\documentclass{amsart}

\usepackage[dvips]{graphics}

\makeatletter
\renewcommand{\p@enumii}{\theenumi-}
\makeatother

\newtheorem{theorem}{Theorem}[section]
\newtheorem{corollary}[theorem]{Corollary}
\newtheorem{lemma}[theorem]{Lemma}
\newtheorem{problem}{Problem}

\theoremstyle{remark}
\newtheorem{remark}[theorem]{Remark}

\numberwithin{equation}{section}

\newcommand{\Diagram}{\mathcal{D}}
\newcommand{\Strip}{\mathcal{S}}
\newcommand{\Twist}{\tau}
\newcommand{\Remainder}{\rho}
\newcommand{\Cancel}{\kappa}
\newcommand{\Edge}{e}
\newcommand{\Edgepath}{\gamma}
\newcommand{\ConstantEdgepaths}{\Gamma_{\mathrm{const}}}
\newcommand{\NonconstantEdgepaths}{\Gamma_{\mathrm{non-const}}}
\newcommand{\EdgepathSystem}{\Gamma}
\newcommand{\BasicEdgepath}{\lambda}
\newcommand{\BasicEdgepathSystem}{\Lambda}

\newcommand{\NumTangles}{N}
\newcommand{\Numer}{P}
\newcommand{\Denom}{Q}
\newcommand{\Slope}{R}

\newcommand{\anglebb}[1]{\langle\langle #1 \rangle\rangle}
\newcommand{\angleb}[1]{\langle #1 \rangle} 
\newcommand{\circlebb}[1]{\langle\langle #1 \rangle\rangle^{\circ}}
\newcommand{\circleb}[1]{\langle #1 \rangle^{\circ}} 

\begin{document}

%
%

\title[Bounds on numerical boundary slopes for Montesinos knots]
{Bounds on numerical boundary slopes\\ for Montesinos knots}

\author{Kazuhiro Ichihara}
\address{
College of General Education, 
Osaka Sangyo University, 
3--1--1 Nakagaito, Daito, Osaka 574--8530, Japan 
}
\email{ichihara@las.osaka-sandai.ac.jp}

\author{%
    Shigeru Mizushima}
\address{%
        Department of Mathematical and Computing Sciences \\
        Tokyo Institute of Technology \\
        12--1 Ohokayama, Meguro \\
        Tokyo 152--8552, Japan}
\email{mizusima@is.titech.ac.jp}

\keywords{boundary slopes, Montesinos knots}
\subjclass[2000]{Primary 57M50; Secondary 57M25}

%
%

\begin{abstract}
We give 
an upper bound on the denominators of numerical boundary slopes and
an upper bound on the differences between two numerical boundary slopes
for Montesinos knots. 
\end{abstract}

\maketitle

%
%

\section{Introduction}

%
%

We consider compact connected surfaces properly embedded 
in compact orientable irreducible $3$-manifolds with single toral boundary, 
which are \textit{essential}, meaning that incompressible and boundary-incompressible. 
The boundary of such a surface consists of a parallel family of 
non-trivial simple closed curves. 
Thus they determine a \textit{slope}, that is, 
the isotopy class of non-trivial simple closed curves. 
This slope is called the \textit{boundary slope} of the surface. 
Boundary slopes of essential surfaces 
have been well-studied, especially, 
in a relation to the study of Dehn surgery on knots. 

Recall that, for the knot exteriors in the $3$-sphere $S^3$, 
the set of slopes is usually identified with the set 
of rational numbers with the infinity $\infty$. 
In fact, such an identification can be done 
by using the standard meridian-longitude system. 
See \cite{Ro} for example. 

In this paper, we study numerical properties of the boundary slopes, 
regarded as rational numbers, for \textit{Montesinos knots}; 
the knots composed by a number of rational tangles. 
Precisely, the aim of this paper is: To give 
(1) an upper bound on the denominator of a boundary slope and
(2) an upper bound on differences between two boundary slopes 
for a Montesinos knot exterior. 
Our bounds are actually described 
in terms of the Euler characteristic and some 
other topological quantity of the surfaces.

For Montesinos knot exteriors, 
Hatcher and Oertel studied 
the boundary slopes in \cite{HO} intensively. 
They gave an algorithm, 
based on the arguments developed originally in \cite{HT}, 
to list all essential surfaces up for a given Montesinos knot exterior. 
Their algorithm has somehow combinatorial workings, 
and in fact, was implemented to a computer program 
by Dunfield described in \cite{Dun}. 
By using this program, we had performed computer-aided experiments, 
and got some observations which suggest the existence of 
such numerical properties for the boundary slopes. 
This is the motivation of our study. 

In the following, 
let $K$ be a Montesinos knot $K(K_1, K_2, \cdots, K_\NumTangles)$,
where $\NumTangles\ge 3$ is the number of tangles and 
each $K_i$ is a non-integral rational number.

%
%

\subsection{Bound on denominator}

We first give an upper bound on the denominators 
of boundary slopes for Montesinos knot exteriors. 

\begin{theorem}
\label{Thm:Denom:UpperBound:Main}
Let $\chi$ be 
the Euler characteristic of the surface corresponding to 
a finite boundary slope $\Slope$ 
and  $\sharp b$ the number of its boundary components. 
Then, 
except for some boundary slopes,
the denominator $\Denom$ of $\Slope$ is bounded as 
\begin{eqnarray}
\Denom \le \frac{-\chi}{\sharp b}.
\label{Eq:Denom:UpperBound:Main-A}
\end{eqnarray}

The exceptions occur from $(-2,3,t)$-pretzel knots for odd $t\ge 3$ or their mirror images.
Some boundary slopes for the knot only satisfies a weaker bound
\begin{eqnarray}
\frac{-\chi}{\sharp b} < \Denom \le \frac{-\chi}{\sharp b} + 1,
\label{Eq:Denom:UpperBound:Main-B}
\end{eqnarray}
though a stronger condition $\sharp b \ge 2$ 
on the number of boundary components holds in these cases. 
\end{theorem}
Here, we remark that, as well as $K(-1/2,1/3,1/t)$, 
for example, Montesinos knots $K((-1/2)+k,(1/3)+l,(1/t)-k-l)$ for $k,l\in \mathbb{Z}$
are also isotopic to the pretzel knot.

From Theorem \ref{Thm:Denom:UpperBound:Main}, 
we have the following corollary immediately. 

\begin{corollary}
\label{Cor:Denom:UpperBound:Main:ByGenus}
Under the same assumption as in Theorem $\ref{Thm:Denom:UpperBound:Main}$, 
and if the surface considered is orientable of genus $g$, 
then the denominator of the boundary slope is bounded as
$\Denom = 1$ if $g = 0 $, 
$\Denom \le 2$ if $g = 1 $, and 
$\Denom \le 2g-1$ if $g \ge 2 $. 
Furthermore, there are no non-torus Montesinos knots 
whose exterior contains essential planar surfaces. 
Thus non-torus Montesinos knots admit no reducible surgery. 
\end{corollary}

The last statement assures that the well-known Cabling Conjecture 
is true for Montesinos knots directly. 
This fact has already been achieved in \cite{EM92} 
as a corollary of the result for strongly invertible knots.



The next corollary is the non-orientable version of the above. 
Recall that a non-orientable surface is called 
of \textit{non-orientable genus} $h$ 
if it contains mutually disjoint $h$ Mobius bands. 

\begin{corollary}
\label{Cor:Denom:UpperBound:Main:ByNonOrientableGenus}
Under the same assumption as in Theorem $\ref{Thm:Denom:UpperBound:Main}$,
and if the surface considered is a non-orientable surface
of non-orientable genus $h$,
then for the denominator $\Denom$ of the boundary slope, we have,
\begin{eqnarray}
\Denom \le \frac{h}{2} +1.
\end{eqnarray}
Moreover if $\sharp b = 1$,
\begin{eqnarray}
\Denom \le h - 1
\end{eqnarray}
holds.
\end{corollary}

%
%

\subsection{Bound on difference}

We next give an upper bound on the ``difference" 
between two boundary slopes for Montesinos knot exteriors 
by a linear function of the ratio $-\chi/\sharp s$ of 
the negative of  the Euler characteristic of the surface 
and the number of sheets. 

The \textit{number of sheets} is the number of pieces of the surface 
in a small neighborhood of a point on a knot. 
If small meridian circles of the knot meet the surface in $m$ points, 
then the number of sheets is $m$.

\begin{theorem}
\label{Thm:Diff:UpperBound:Main}
Let $\chi_i$ be 
the Euler characteristic of the surface 
corresponding to a finite boundary slope $\Slope_i$ 
and $\sharp s_i$ its number of sheets, for $i=1,2$ respectively. 
Then
the difference $|\Slope_1-\Slope_2|$ between the boundary slopes $\Slope_1$ and $\Slope_2$
is bounded as
\begin{eqnarray}
|\Slope_1-\Slope_2|\le 2\,(\frac{-\chi_1}{\sharp s_1}+\frac{-\chi_2}{\sharp s_2})+4
.
\label{Eq:Diff:UpperBound:Main}
\end{eqnarray}
\end{theorem}

This inequality (\ref{Eq:Diff:UpperBound:Main}) can be rewritten as 
\begin{eqnarray}
\Delta(\Slope_1,\Slope_2)\le 2\,(\Denom_2 \frac{-\chi_1}{\sharp b_1}+\Denom_1 \frac{-\chi_2}{\sharp b_2})+4 \Denom_1 \Denom_2
,
\label{Eq:Dist:UpperBound:Main-C}
\end{eqnarray}
which may be preferable for understanding the meaning 
from the geometric viewpoint. 
Here $\Delta(\Slope_1,\Slope_2)$ denotes the \textit{distance} 
between the slopes $\Slope_1$ and $\Slope_2$, 
which is defined to be the minimal geometric intersection number 
of the simple closed curves representing $\Slope_1$ and $\Slope_2$. 
Recall that if $\Slope_i$ is expressed by 
an irreducible fraction $\Numer_i/\Denom_i$ for $i=1,2$, 
then $\Delta(\Slope_1,\Slope_2)$ is equal to $| \Numer_1 \Denom_2 - \Numer_2 \Denom_1 |$. 
However, in the algorithm of Hatcher and Oertel, 
$\Slope_i=\Numer_i/\Denom_i$ and $\sharp s_i$ 
play significant roles rather than $\Numer_i$, $\Denom_i$ and $\sharp b_i$.
Hence, in the light of the algorithm, it seems natural to consider 
the difference $|\Slope_1-\Slope_2|$ and $-\chi_i/\sharp s_i$. 
Note that, in particular case that both $\Slope_1$ and $\Slope_2$ are integers, 
$|\Slope_1-\Slope_2|$ coincides with $\Delta(\Slope_1,\Slope_2)$, 
and we have an upper bound of the distance simultaneously.

From Theorem $\ref{Thm:Diff:UpperBound:Main}$, 
we have three corollaries as follows. 

When the surface are both orientable, we immediately have the following. 

\begin{corollary}
\label{Cor:Diff:UpperboundByGenus}
Under the same assumption as in Theorem $\ref{Thm:Diff:UpperBound:Main}$,
and if the surfaces considered are both orientable surfaces
of genera $g_1$ and $g_2$ respectively, then,
the difference $|\Slope_1-\Slope_2|$ between the boundary slopes $\Slope_1$ and $\Slope_2$
is bounded as
\begin{eqnarray}
|\Slope_1-\Slope_2| \le 4 \,( g_1 + g_2 ) .
\label{Eq:Diff:UpperBound:Main-B:ByGenus}
\end{eqnarray}
\end{corollary}

With respect to a linear bound on the difference, 
or a somehow irregular quadratic bound on the distance by Euler characteristics,
the following corollary is easily obtained from Theorem \ref{Thm:Diff:UpperBound:Main}.
Though the bounds may not be sharp for Montesinos knots with $\NumTangles\ge 3$,
the equality holds for boundary slopes of the trefoil knot.


\begin{corollary}
\label{Cor:DiffDist:UpperBound}
For two boundary slopes and their corresponding essential surfaces,
we have the inequality
\begin{eqnarray}
|\Slope_1-\Slope_2|\le 6\,(\frac{-\chi_1}{\sharp s_1}+\frac{-\chi_2}{\sharp s_2}).
\label{Eq:Diff:UpperBound:Linear}
\end{eqnarray}
This is equivalent to the inequality
\begin{eqnarray}
\Delta(\Slope_1,\Slope_2)\le 6\,(\Denom_2 \frac{-\chi_1}{\sharp b_1}+\Denom_1 \frac{-\chi_2}{\sharp b_2}).
\label{Eq:Dist:UpperBound:Semi-Linear}
\end{eqnarray}
\end{corollary}

Regarding the upper bound of the distance or difference by the product of Euler characteristics, we have the following.
Though the bounds may not be sharp for Montesinos knots with $\NumTangles\ge3$ tangles,
the equality holds for boundary slopes of the figure eight knot.


\begin{corollary}
\label{Cor:Main3b}
If both of the Euler characteristics are negative, then we have
\begin{eqnarray}
\Delta(\Slope_1,\Slope_2)\le 8\cdot \frac{-\chi_1}{\sharp b_1}\cdot \frac{-\chi_2}{\sharp b_2}.
\label{Eq:Dist:UpperBound:Quadratic}
\end{eqnarray}
This is equivalent to the inequality
\begin{eqnarray}
|\Slope_1-\Slope_2|\le 8\cdot \frac{-\chi_1}{\sharp s_1}\cdot \frac{-\chi_2}{\sharp s_2}.
\label{Eq:Diff:UpperBound:Quadratic}
\end{eqnarray}
\end{corollary}

%
%

This paper is organized as follows. 
We review the algorithm in \cite{HO} in Section 2 and 
prepare some formulae in Section 3. 
Then, Section 4 and 5 are devoted to giving proofs of 
Theorem \ref{Thm:Denom:UpperBound:Main} and 
\ref{Thm:Diff:UpperBound:Main} respectively.
In the last section, 
a brief review on related known results is given, 
and some open problems are stated.

%
%

\section*{Acknowledgments}
The authors would like to thank Professor Sadayoshi Kojima 
for his helpful suggestions about earlier drafts. 
They also thank to Professor Masakazu Teragaito for 
letting them know the related paper \cite{Te}. 

%
%

\section{Algorithm of Hatcher and Oertel}

The proofs of both Theorem \ref{Thm:Denom:UpperBound:Main} and \ref{Thm:Diff:UpperBound:Main} deeply depend on the algorithm in \cite{HO}.
Hence, in this section,
we review its workings of enumerating all boundary slopes.


\subsection*{Montesinos knot}

As mentioned in the introduction,
we assume that the knot $K$ is a Montesinos knot $K(K_1,K_2,\ldots,K_\NumTangles)$, where each $K_i$ is a non-integral fraction and $\NumTangles\ge 3$ in this article.
By the assumption, 
we normalize Montesinos knots and
eliminate two-bridge knots from the argument.
This is because boundary slopes are enumerated for two-bridge knots in \cite{HT},
and two-bridge knot case is excluded in \cite{HO}.
Results for two-bridge knots similar to our main results are obtained by \cite{HT}.

Since a knot in this article is basically a Montesinos knot,
we use the term ``tangle'' as a rational tangle if not mentioned otherwise particularly.


\subsection*{Decomposition}

First, we regard $S^3$ including a Montesinos knot $K$ as the union of a $\NumTangles$-tuple of 3-balls $B_i~(i=1, 2, \ldots, \NumTangles)$ in $S^3$
with following properties.
The interiors of all $B_i$'s are disjoint.
The intersection of all boundaries $\partial B_i$'s is a circle, which is called the {\em axis} of the knot $K$.
Each $\partial B_i$ is divided into two hemispheres by the axis,
and the right hemisphere of $\partial B_i$ coincides with the left hemisphere of $\partial B_{i+1}$ (indices are taken modulo $\NumTangles$).
Each ball $B_i$ includes the rational tangle $K_i$ of the Montesinos knot $K$.

By this decomposition,
a properly embedded essential surface $F$ is also decomposed into a $\NumTangles$-tuple of subsurfaces $S_i$ in $B_i$.
The boundary of $S_i$ is the union of the tangle $K_i$ in the interior of $B_i$ and a curve system on the four-punctured sphere $\partial B_i\setminus K_i$,
where a {\em curve system} means the union of disjoint circles and arcs connecting distinct punctures.
A simple example of a curve system is a {\em $p/q$-tangle} drawn on a four-punctured sphere. 
It is denoted by $\anglebb{p/q}$.
Note that a tangle usually means two strings in a 3-ball with their four ends fixed on the boundary sphere,
in some cases,
we use ``tangle'' as a rational tangle projected to, or drawn on the boundary disjointly.
Another example of a curve system is a {\em $p/q$-circle},
which is a non-trivial circle disjoint from $p/q$-tangle in a sphere,
and is denoted by $\circlebb{p/q}$.


\subsection*{Subsurfaces}

In the argument,
subsurfaces $S_i$'s are arranged to sit in a standard position by isotopy,
and are restricted to be saddle subsurfaces or cap subsurfaces
as shown in \cite{HO}.

The simplest example of a subsurface $S_i$ is the direct product set of the curve system $\anglebb{p/q}$ in a four-punctured level sphere with an interval.
Also this subsurface is denoted by $\anglebb{p/q}$.
This subsurface $S_i$ is topologically the union of two disks,
and is called {\em base disks}
since every subsurface can be regarded to include these kind of disks.

\begin{figure}[hbt]
 \begin{picture}(330,130)
 \put(40,20){\scalebox{1.0}{\includegraphics{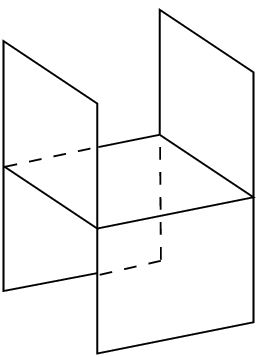}}} 
 \put(0,0){(a) A saddle subsurface $\anglebb{\infty}$--$\anglebb{0}$}
 \put(220,20){\scalebox{0.8}{\includegraphics{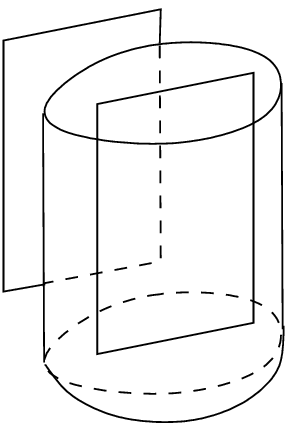}}} 
 \put(180,0){(b) A cap subsurface $\anglebb{0}+\circlebb{0}$}
 \end{picture}
 \caption{Examples of subsurfaces}
 \label{Fig:Subsurfaces}
\end{figure}

An example of a saddle subsurface $S_i$ is constructed by connecting two surfaces by a saddle,
where both of the two surfaces are ``base disks'' $\anglebb{p/q}$ and $\anglebb{r/s}$ described above
for $p/q$ and $r/s$ satisfying $|ps-qr|=1$.
See Figure \ref{Fig:Subsurfaces}(a).
A saddle is a disk on a level sphere
bounded by a simple closed curve 
made of $p/q$-tangle, $r/s$-tangle and four punctures.
Though there are two choices of disks bounded by the simple closed curve,
the choice does not matter in our later argument.
Note that $|ps-qr|=1$ ensures that $p/q$-tangle and $r/s$-tangle are disjoint in a level sphere.
This subsurface is denoted by $\anglebb{p/q}$\,--\,$\anglebb{r/s}$.

We can construct the disjoint union of 
$k$ parallel copies of a subsurface $\anglebb{r/s}$
and $l$ parallel copies of a saddle subsurface $\anglebb{p/q}$\,--\,$\anglebb{r/s}$.
This subsurface $S_i$ is denoted by $(k\anglebb{p/q}+l\anglebb{r/s})$\,--\,$(k+l)\anglebb{r/s}$.

For a sequence
$(k\anglebb{p_{j}/q_{j}}+l\anglebb{p_{j-1}/q_{j-1}})$\,--\,$(k+l)\anglebb{p_{j-1}/q_{j-1}}$\,--\,$\ldots$\,--\,$(k+l)\anglebb{p_{1}/q_{1}}$,
we can construct a subsurface by preparing components corresponding to each pair of successive two curve systems in the sequence and gluing them together according to the sequence.
The rightmost curve system $p_{1}/q_{1}$ is required to coincide with $K_i$ so that the boundary of $S_i$ includes the rational tangle $K_i$.
We regard $\anglebb{p_{1}/q_{1}}$ as the starting point of the sequence.
Furthermore, we describe curve systems from right to left in a sequence of curve systems as above.


A cap subsurface is constructed as follows.
We prepare a curve system consisting of $k$ parallel copies of $p/q$-tangle and $l$ parallel copies of $p/q$-circles.
Next, we take a direct product of the curve system with the interval,
and arrange the product to lie inside the ball $B_i$
so that one of the two boundary level spheres is placed at $\partial B_i$.
Then, $p/q$-circles of the inner boundary of the product are capped by disks.
See Figure \ref{Fig:Subsurfaces}(b).
Both this curve system and the cap subsurface are denoted by $k\anglebb{p/q}+l\circlebb{p/q}$.

For any type of subsurface, the leftmost curve system of a sequence represents the curve system $S_i \cap \partial B_i$.

Here, we note that some surfaces may correspond to the same representation by $\anglebb{}$ and $\circlebb{}$,
because of the two choices of saddles described above.
Though, 
the ambiguity does not cause trouble,
and we regard a representation as if it corresponded to a surface. 

\subsection*{$abc$-coordinates}

\begin{figure}[htb]
   \begin{center}
    \begin{picture}(90,105)
     \put(0,0){\scalebox{0.5}{\includegraphics{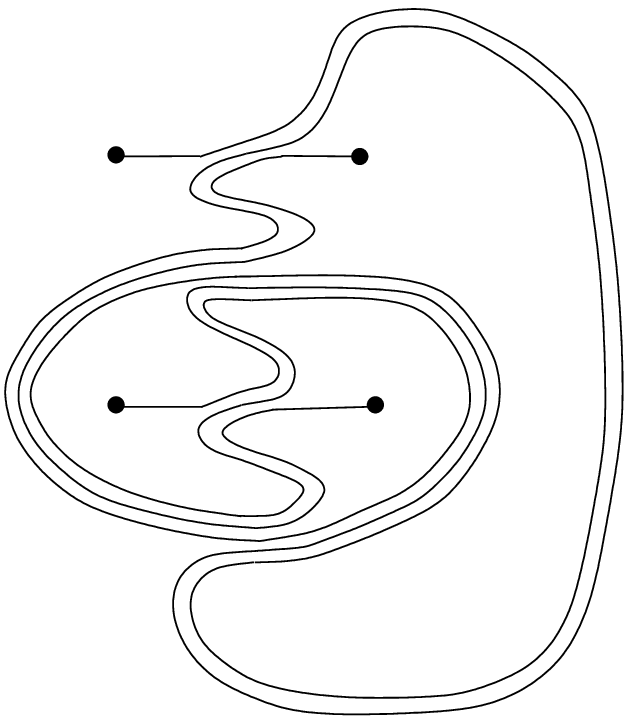}}} 
     \put(18,48){$a$}
     \put(50,48){$a$}
     \put(18,85){$a$}
     \put(50,85){$a$}
     \put(66,25){$b$}
     \put(91,30){$c$}
    \end{picture}
   \end{center}
   \caption{A curve system with $(a,b,c)=(1,3,2)$}
   \label{Fig:CurveSystem}
  \end{figure}

As illustrated in Figure \ref{Fig:CurveSystem},
we normalize a curve system to a standard form.
The curve system in the standard form is represented by integers $a$, $b$ and $c$.
$a$, $b$ and $c$ denote the number of subarcs of the curve system lying in a particular region as in the figure.
The number of subarcs around the axis, denoted by $c$, can be negative. 
If $c$ is negative, $(a,b,c)$ represents a curve system which is the mirror image of the curve system represented by $(a,b,|c|)$ taken with respect to the axis.
Thus, $(a,b,c)$ represents a curve system and these are called {\em $abc$-coordinates} of a curve system.

For example, the coordinates of 
a $p/q$-tangle and a $p/q$-circle are $(1,q-1,p)$ and $(0,q,p)$ respectively.
The coordinates of the disjoint union of two curve systems
are calculated as the vector sum of the coordinates of both curve systems.

Though there is no explicit description about $\anglebb{\infty}$ in \cite{HO},
we only have to introduce the $d$-coordinate for the number of $\infty$-tangles.

With the $abc$-coordinates,
we can explicitly describe conditions for subsurfaces $S_i$'s to be glued consistently.
For the $abc$-coordinates of the last curve systems for all subsurfaces,
one of the conditions is that $a$-coordinates are the same for all $S_i$'s and so are $b$-coordinates.
The other condition is that $c$-coordinates for all $S_i$'s sum up to exactly $0$.



\subsection*{$uv$-coordinates}

By projectifying the $abc$-coordinates to the {\em $uv$-coordinates}
by $u=b/(a+b)$ and $v=c/(a+b)$,
we can make the subsequent argument simpler.
An important fact is that a curve system consisting of parallel copies of a curve system of arbitrary multiplicity is mapped to the same point in the $uv$-plane.

A curve system $\anglebb{p/q}$ is projected to $(u,v)=((q-1)/q, p/q))$,
which is denoted by $\angleb{p/q}$.
%
%
A curve system $(k\anglebb{p/q}+l\anglebb{r/s})$ has coordinates
$(a,b,c)=(k+l,k(q-1)+l(s-1),kp+lr)$,
and thus is projected to the point $(u,v)=(kq/(kq+ls))((q-1)/q,p/q)+(ls/(kq+ls))((s-1)/s,r/s)$.
Thus, for fixed $p$, $q$, $r$ and $s$, the $uv$-coordinates depend on only the ratio of $k$ to $l$.
After projectification,
the curve system is represented as a point on the segment connecting $\angleb{p/q}$ and $\angleb{r/s}$.
This point is denoted by
$((k/(k+l))\angleb{p/q}+(l/(k+l))\angleb{r/s})$.
Note that we use the ratio $k/(k+l)$ instead of $kq/(kq+ls)$,
since it is suitable in later calculation.

A curve system $\circlebb{p/q}$ is projected to $(u,v)=(1,p/q)$ on a vertical line $u=1$, and is denoted by $\circleb{p/q}$.
%
A curve system $(k\anglebb{p/q}+l\circlebb{p/q})$
has coordinates $(a,b,c)=(k,k(q-1)+lq,kp+lp)$,
and thus is projected to $(u,v)=(k/(k+l))\cdot((q-1)/q,p/q)+(l/(k+l))\cdot(1,p/q)$
on the horizontal segment connecting $\angleb{p/q}$ and $\circleb{p/q}$.
Note that, for fixed $p$ and $q$,
the $uv$-coordinates depend on only the ratio of $k$ to $l$.
This point is denoted by $((k/(k+l))\angleb{p/q}+(l/(k+l))\circleb{p/q})$.

After projectification,
from a sequence of points in the $abc$-space representing a subsurface $S_i$,
we obtain a sequence of points on the $uv$-plane for the subsurface.

The conditions for subsurfaces to be glued together consistently
with respect to the $abc$-coordinates
are translated into the conditions of the $uv$-coordinates.
For the $uv$-coordinates of the last one of the sequence of curve systems for all subsurfaces,
$u$-coordinates are the same for all subsurfaces
and $v$-coordinates for all subsurfaces sum up to $0$.


\subsection*{The diagram and edgepaths}

A subsurface $S_i$ is formally related to 
a piecewise linear path in the $uv$-plane.
Such a path is called an {\em edgepath}, and
we will often use $\Edgepath$ as the symbol for it.
Edgepaths lie on a ``diagram'' described as follows.

The {\em diagram $\Diagram$} is a graph on the $u$-$v$ plane.
A vertex is 
a point $\angleb{p/q}$, whose coordinates are $(u,v)=((q-1)/q,p/q)$,
a point $\circleb{p/q}$, whose coordinates are $(u,v)=(1,p/q)$, where $p/q$ is an irreducible fraction, 
or a point $\angleb{\infty}$, whose coordinates are $(u,v)=(-1,0)$.
If two vertices $\angleb{p/q}$ and $\angleb{r/s}$ satisfy the condition $|ps-qr|=1$, the two vertices are connected by a segment.
This segment is one of the two types of edges of the diagram and is denoted by $\angleb{p/q}$\,--\,$\angleb{r/s}$.
We call the segment a {\em non-horizontal edge}.
In particular, for an integer $z$, 
there are edges $\angleb{z}$\,--\,$\angleb{z+1}$ and $\angleb{\infty}$\,--\,$\angleb{z}$.
The former is called a {\em vertical edge} since it is a segment of the vertical line $u=0$. The latter is called an {\em $\infty$-edge}.
%
The other type of edge is called a {\em horizontal edge},
which connects vertices $\angleb{p/q}=((q-1)/q,p/q)$ and $\circleb{p/q}=(1,p/q)$.
This edge is denoted by $\angleb{p/q}$\,--\,$\circleb{p/q}$.
Note that though the edge $\angleb{\infty}$\,--\,$\angleb{0}$ is horizontal in the usual sense,
we regard the edge as a non-horizontal edge rather than a horizontal edge
for ease in our later argument.
The region $-1\le u \le 1$ is triangulated by these kinds of edges,
though the triangulation is not locally finite.
In particular, the part of the diagram lying in the strip $0\le u\le 1$ 
is denoted by $\Strip$.

  \begin{figure}[htb]
   \begin{center}
    \begin{picture}(70,150)
     \put(0,0){\scalebox{0.7}{\includegraphics{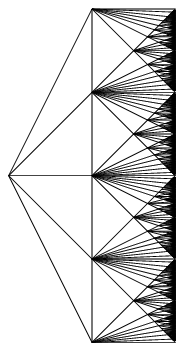}}}
     \put(49,135){\rotatebox{90}{\scalebox{1.0}{$\cdots$}}}
     \put(49,-10){\rotatebox{90}{\scalebox{1.0}{$\cdots$}}}
     \put(42,70){\scalebox{1.5}{\vector(1,0){30}}}
     \put(42,70){\scalebox{1.5}{\vector(0,1){50}}}
     \put(92,68){$u$}
     \put(40,148){$v$}
     \put(33,73){$O$}
    \end{picture}
   \end{center}
   \caption{The diagram $\Diagram$}
   \label{Fig:Diagram}
  \end{figure}

The edgepath $\Edgepath_i$ of a cap subsurface $S_i$ is a point on the horizontal segment $\angleb{K_i}$\,--\,$\circleb{K_i}$.
The edgepath $\Edgepath_i$ of a saddle subsurface $S_i$ is a piecewise linear path starting from the vertex $\angleb{K_i}$.
The endpoint of an edgepath is either of a vertex of the diagram or 
a point on an edge of the diagram.
Hence, the last edge of an edgepath may be a part of a non-horizontal edge.
We call such an edge a {\em partial edge}.
In comparison with this, we use the term {\em complete edge} to express the whole of a non-horizontal edge.
An edgepath consisting of only one point is called a {\em constant edgepath}.
The other type of edgepath is called a {\em non-constant edgepath}.


\subsection*{Edgepath systems}

By collecting edgepaths for subsurfaces $S_i$,
we can represent the original surface $F$.
We call this kind of $\NumTangles$-tuple ($\Edgepath_1$, $\Edgepath_2$, $\ldots$, $\Edgepath_\NumTangles$) an {\em edgepath system}
and will often use $\EdgepathSystem$ as the symbol for it.

Conversely,
for an appropriate edgepath system,
by unprojecting all vertices in its edgepaths in the $uv$-plane
to integral points in the $abc$-space with the common $a$-coordinate,
we can construct subsurfaces and a surface $F$, though some ambiguity remains.


The set of the edgepath systems are divided into three classes 
according to the common $u$-coordinate of the endpoints of edgepaths in their edgepath system.
An edgepath system and the corresponding surface are said to be {\em type I}, {\em type II} or {\em type III},
if all edgepaths in the edgepath system end at $u>0$, $u=0$ or $u<0$ respectively.


\subsection*{Candidate surfaces}

In the enumeration of boundary slopes, 
we first list candidates for essential surfaces,
and then omit compressible surfaces from the candidates.
Precise conditions for an edgepath system $\EdgepathSystem=(\Edgepath_1, \Edgepath_2,\ldots ,\Edgepath_\NumTangles)$ to be an edgepath system of a {\em candidate surface} are given in \cite{HO} as follows.
\begin{itemize}
\item[(E1)] The starting point of $\Edgepath_i$ lies on the edge $\angleb{K_i}$\,--\,$\circleb{K_i}$, and if this starting point is not the vertex $\angleb{K_i}$, then the edgepath $\Edgepath_i$ is constant.
\item[(E2)] $\Edgepath_i$ is minimal, i.e., it never stops and retraces itself, nor does it ever go along two sides of the same triangle of $\Diagram$ in succession.
\item[(E3)] The ending points of the $\Edgepath_i$'s are rational points of $\Diagram$ which all lie on one vertical line and whose vertical coordinates add up to zero.
\item[(E4)] $\Edgepath_i$ proceeds monotonically from right to left, ``monotonically'' in the weak sense that motion along vertical edges is permitted. 
\end{itemize}


\subsection*{Basic edgepath systems and consistency on gluing}

The enumeration is performed by use of basic edgepath systems.
A {\em basic edgepath} is an edgepath which
starts from a vertex $\angleb{p/q}$,
goes leftwards monotonically,
and ends at the time when the edgepath first reaches $u=0$.
Moreover, a {\em basic edgepath system}
is an edgepath system which consists of basic edgepaths.
We will often use $\BasicEdgepath$ and $\BasicEdgepathSystem$ as the symbols for
a basic edgepath and a basic edgepath system.

In order to seek type I edgepaths,
it is convenient to introduce an {\em extended basic edgepath} $\widetilde{\BasicEdgepath}$
which is obtained 
by connecting a horizontal segment $\angleb{p/q}$\,--\,$\circleb{p/q}$
to the starting point $\angleb{p/q}$ of a basic edgepath $\BasicEdgepath$.
We define an {\em extended basic edgepath system} $\widetilde{\BasicEdgepathSystem}$ similarly.

Sometimes, we regard 
an edgepath $\Edgepath$ as a function from an interval in $\mathbb{R}$ to $\mathbb{R}$, 
which maps $u$-coordinate to $v$-coordinate,
and then allow ourselves to use expressions like $\Edgepath(u)$,
although we cannot define its value for $u=0$ if the edgepath includes vertical edges.
The function is piecewise linear.
Similar notation is used for
a basic edgepath and an extended basic edgepath.
With the notation, a condition of consistency on gluing in (E3) can be described as 
\begin{eqnarray}
\label{Eq:EquationForEdgepathSystem}
\sum_{i=1}^{\NumTangles} \Edgepath_i(u)=0.
\end{eqnarray}
Especially for type I surfaces,
we need to solve (\ref{Eq:EquationForEdgepathSystem}) for some extended basic edgepath system
to determine the common $u$-coordinate of the endpoints of its edgepath system.


\subsection*{Enumeration}
We have finished introducing notions used in the algorithm in \cite{HO}.
Now, we review its workings.
All boundary slopes are enumerated as follows.

All basic edgepath systems for the Montesinos knot $K$ are enumerated first.
Then, type I, type II and type III candidate edgepath systems $\EdgepathSystem$ are obtained for each basic edgepath system $\BasicEdgepathSystem$ .
A type I edgepath system is obtained by solving the equation (\ref{Eq:EquationForEdgepathSystem}) for the extended basic edgepath system $\widetilde{\BasicEdgepathSystem}$.
For a solution $u_0$ of the equation, we construct an edgepath system as follows.
Let $\widetilde{\BasicEdgepath_i}$ be the $i$-th extended basic edgepath
starting from $\angleb{K_i}$ where $K_i=p_i/q_i$.
If $(q_i-1)/q_i < u_0$, then the line $u=u_0$ intersects with the horizontal edge of the extended basic edgepath.
Therefore, we prepare a constant edgepath with a single point $(u,v)=(u_0,K_i)$.
If otherwise, the line $u=u_0$ intersects with the non-horizontal part of the extended basic edgepath.
Hence, we cut out an edgepath $\Edgepath_i$ 
starting from $\angleb{K_i}$ and ending at $u=u_0$,
from the original basic edgepath $\BasicEdgepath_i$.
The edgepath system is obtained by collecting all such edgepaths.
A type II edgepath system $\EdgepathSystem$ is obtained by adding vertical edges to the basic edgepaths of the basic edgepath system $\BasicEdgepathSystem$
so that $v$ coordinates of endpoints of the edgepath system $\EdgepathSystem$ sum up to $0$.
%
A type III edgepath system $\EdgepathSystem$ is obtained by adding an $\infty$-edge to each basic edgepath $\BasicEdgepath_i$.
Thus, we can enumerate all candidate surfaces.
The process is completed in finite time.

After enumerating all the candidate edgepath systems,
we verify their incompressibility.
Detailed conditions for the edgepath system of a candidate surface
to be incompressible are also described in \cite{HO}.
By the conditions, we can eliminate compressible surfaces from the set of candidate surfaces, and complete the enumeration of essential surfaces.
Though, the conditions are not so crucial in this paper
and we hardly make use of the conditions.
Besides, the determination of the orientability is omitted in \cite{HO}.
It must be performed by oneself if necessary.


%
%

\section{Preparation}
\label{Sec:Preparation}

In this section,
we prepare some formulae for concrete calculation in the subsequent sections.
We also introduce an operation named ``simplification''.


\subsection*{Notation}

We first give some notation about edgepaths and edgepath systems.
For an edgepath $\Edgepath$,
symbols $\Edgepath_{>0}$ and $\Edgepath_{\ge 0}$ denote
a part of $\Edgepath$ inside the region $u>0$ and $u\ge 0$ respectively.
A part of $\Edgepath$ consisting of vertical edges
is denoted by $\Edgepath_{=0}$.

For type II and type III edgepath systems
let $\Edgepath(+0)$ denote the value of $v$ at the moment when $u$-coordinate reaches $0$.
Furthermore, for an edgepath system $\EdgepathSystem=(\Edgepath_1, \Edgepath_2,\ldots ,\Edgepath_\NumTangles)$, let $\EdgepathSystem(+0)$ denote the sum $\sum_{i=1}^{\NumTangles} \Edgepath_i(+0)$.


\subsection*{Signs of edges}
A complete non-horizontal edge $\Edge=\angleb{p/q}$\,--\,$\angleb{r/s}$
is said to be {\em increasing} or {\em decreasing}
if the $v$-coordinate of a point increases or decreases respectively when it goes from $\angleb{r/s}$ to $\angleb{p/q}$ along $\Edge$.

We define the sign of the edge $\Edge$ to be $+1$ or $-1$
according to whether the edge is increasing or decreasing respectively.
The sign of the edge $\Edge$ is denoted by $\sigma(\Edge)$
and is calculated by $(ps-qr)$.
The sign of a partial edge is defined as the sign of the complete edge including the partial edge.


\subsection*{Lengths of edgepaths}

Next, we define the length of an edgepath.
The lengths of a single point and a complete edge are defined to be $0$ and $1$ respectively.
The length of a partial edge $\Edge=(k/(k+l) \angleb{p/q}+l/(k+l) \angleb{r/s})$\,--\,$\angleb{r/s}$ is $k/(k+l)$.
Note that the ratio of the Euclidean length on the $uv$-plane of the partial edge 
to that of the complete edge $\angleb{p/q}$\,--\,$\angleb{r/s}$ is $kq/(kq+ls)$ as calculated in the previous section.
Thus, the length of a partial edge does not coincide with the ratio generally.
The length of an edge $\Edge$ is denoted by $|\Edge|$.
The length of an edgepath is the sum of the lengths of the edges in the edgepath.
A constant edgepath is of length $0$.
The length of an edgepath $\Edgepath$ is denoted by $|\Edgepath|$.

Here, we prepare another formula of the length of a partial edge.
For a partial edge $e$ of a complete edge $\angleb{p/q}$\,--\,$\angleb{r/s}$,
assume that $u$-coordinate of the endpoint of the partial edge is $u_0$.
If a curve system $k\anglebb{p/q}+l\anglebb{r/s}$ corresponds to the endpoint,
we have $abc$-coordinates $(a_0,b_0,c_0)=k(1,q-1,p)+l(1,s-1,r)$ and $u$-coordinate $u_0=b_0/(a_0+b_0)=\{k(q-1)+l(s-1)\}/(kq+ls)$,
and hence,
\begin{eqnarray}
|e|&=& \frac{k}{k+l}=\frac{1+s(u_0-1)}{(s-q)(u_0-1)}
\label{Eq:Formula:LengthOfPartialEdge}
.
\end{eqnarray}


\subsection*{Boundary slopes and twists of surfaces}

The boundary slope of a surface is calculated via the total number of twists,
which we call {\em twist} for short.
Roughly, the twist $\Twist(F)$ of a surface $F$ is a variation of the numerical boundary slope, which fits with the algorithm.
With the twist, the boundary slope $\Slope$ is calculated by $\Slope = \Twist(F) - \Twist(F_S)$
where $F_S$ is a Seifert surface in the list of candidate surfaces of the knot $K$.

\begin{figure}[htb]
   \begin{center}
    \begin{picture}(70,102)
     \put(0,0){\scalebox{0.6}{\includegraphics{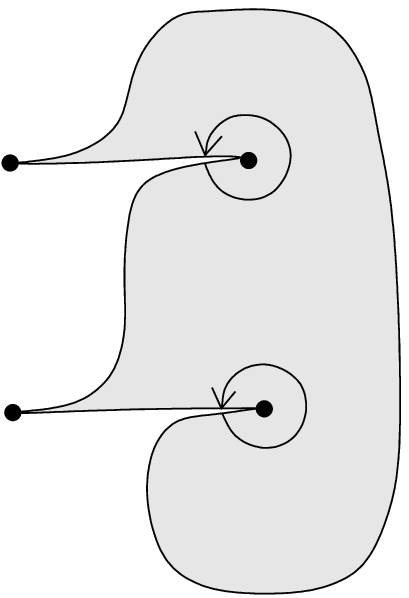}}} 
    \end{picture}
   \end{center}
   \caption{A saddle corresponding to $\angleb{0}$\,--\,$\angleb{1}$}
   \label{Fig:Twist}
  \end{figure}

We define the twist of a subsurface first.
Base disks have twist $0$.
For a non-$\infty$-edge $\anglebb{p/q}$\,--\,$\anglebb{r/s}$,
if we draw both tangles in standard position as in Figure \ref{Fig:CurveSystem},
the saddle for the edge surrounds two of the four punctures
(see Figure \ref{Fig:Twist}).
This means that
two of the four boundary arcs of a saddle component of a subsurface revolve once around the strands of the tangle.
Hence, the edge of the edgepath contributes $\pm 2$ to the twist.
The sign of the value is determined by whether the boundary arc revolves in clockwise or counter-clockwise direction,
and the sign coincides with $-\sigma(\Edge)$ for an edge $\Edge$.
For a partial edge $(k\anglebb{p/q}+l\anglebb{r/s})$\,--\,$(k+l)\anglebb{r/s}$,
at two of the four boundary arcs of a component of a subsurface,
$k$ of $(k+l)$ sheets go around the strands once.
Hence, the partial edge contributes $\pm 2 k/(k+l)$ to the twist.
Note that $k/(k+l)$ coincides with the length of the partial edge.
The twist of a subsurface is the sum of the twists of the edges of the edgepath corresponding to the subsurface.
Naturally, the twist of a cap subsurface is $0$, since the corresponding edgepath is a point.

The twist of the surface $F$ is the sum of the twists of its subsurfaces $S_i$.
The precise definition of the twist is
\begin{eqnarray}
\Twist(F)
&=&
\sum_{i=1}^{\NumTangles}
~
\left\{
 \begin{array}{l}
  0 \\
  ~~~~~(\textrm{  if $\Edgepath_i$ is constant }) \\
  \sum_{\Edge_{i,j}\in \Edgepath_i}
  \left\{
   \begin{array}{l}
    0 \\
    ~~~~~(\textrm{  if $\Edge_{i,j}$ is an $\infty$-edge.} ) \\
    -2\,\sigma(\Edge_{i,j})~|\Edge_{i,j}| \\
    ~~~~~(\textrm{  otherwise. } )\\
   \end{array}
  \right. \\
  ~~~~~(\textrm{  if $\Edgepath_i$ is non-constant }) \\  
 \end{array}
\right.
.
\label{Eq:Formula:Twist}
\end{eqnarray}
Though the twist is originally defined for a surface,
it is well-defined for an edgepath system.


\subsection*{Surfaces with the same boundary slope}

We think about surfaces with the same boundary slope.
In the proofs of the theorems,
only a surface of minimal $-\chi/\sharp s$ is important
among such surfaces sharing the common boundary slope.
Therefore, we take the surface of minimal $-\chi/\sharp s$ as their representative,
ignore the others and will simplify the subsequent argument,
especially in Section \ref{Sec:ABoundOnTheDifference}.
We call this operation {\em simplification}.

By the simplification, we ignore (1) most of type I surfaces corresponding to non-isolated solutions of (\ref{Eq:EquationForEdgepathSystem}), (2) type II surfaces with redundant vertical edges, (3) type III surfaces with partial $\infty$-edges and (4) augmented type III surfaces mentioned in \cite{HO}.

\subsection*{Euler characteristics}

Instead of the Euler characteristic itself,
formulae for calculating $-\chi/\sharp s$ are given
since they are more suitable.
Note that 
though the Euler characteristic itself is not well-defined for an edgepath system,
so is the quantity $-\chi/\sharp s$.


To construct a type III surface,
we have base disks for each tangle first,
add saddles according to non-$\infty$-edges in the edgepath system,
add also saddles according to $\infty$-edges,
and then glue $S_i$'s together at arcs which are halves of $\infty$-tangle on $\partial B_i$.
The Euler characteristic of the surface $F$ so obtained is calculated by
$\chi=
\sum_{i=1}^{\NumTangles} (2\cdot \sharp s)
-\sum_{i=1}^{\NumTangles} (|\Edgepath_{i,>0}|\cdot \sharp s)
-\sum_{i=1}^{\NumTangles} (1\cdot \sharp s)
-\NumTangles \cdot \sharp s
$.
Thus,
\begin{eqnarray}
\frac{-\chi}{\sharp s}
&=&
 \sum_{i=1}^{\NumTangles}|\Edgepath_{i,>0} |
\label{Eq:Formula:EulerCharTypeIII}
.
\end{eqnarray}

\bigskip


To construct a type II surface,
we have base disks first, 
add saddles for the basic edgepath,
add also saddles for vertical edges,
and then glue $S_i$'s together at integral tangles on $S_i \cap \partial B_i$.
Euler characteristic is
$\chi=
\sum_{i=1}^{\NumTangles} (2\cdot \sharp s)
-\sum_{i=1}^{\NumTangles} (|\Edgepath_{i,>0}|\cdot \sharp s)
-\sum_{i=1}^{\NumTangles} (|\Edgepath_{i,=0}|\cdot \sharp s)
-2(\NumTangles-1) \cdot \sharp s
$.
Thus,
\begin{eqnarray}
\frac{-\chi}{\sharp s}
&=&
 \sum_{i=1}^{\NumTangles}(|\Edgepath_{i,>0} |)
 +|\EdgepathSystem(+0) |-2
\label{Eq:Formula:EulerCharTypeII}
.
\end{eqnarray}

\bigskip


To construct a type I surface,
we have base disks first,
add caps for constant edgepaths,
add saddles for non-constant edgepaths,
and then glue $S_i$'s at curve systems on $\partial B_i$'s.

Assume that $C$ is a component of a subsurface described by $(k\anglebb{p/q}+l\anglebb{r/s})$\,--\,$(k+l)\anglebb{r/s}$.
$k$ saddles are included in the component $C$,
and contribute $-k$ to the Euler characteristic.
Since $\sharp s=k+l$,
the contribution to $-\chi/\sharp s$ by the partial edge $(k/(k+l)\angleb{p/q}+l/(k+l)\angleb{r/s})$\,--\,$\angleb{r/s}$ is $k/(k+l)$,
which coincides with the length of the partial edge $e$ by definition.

Next assume that $C$ is a cap subsurface described by $k\anglebb{p/q}+l\circlebb{p/q}$.
$l$ caps are included in the component $C$,
and contribute $+l$ to Euler characteristic.
By $\sharp s=k$,
the contribution to $-\chi/\sharp s$ by the constant edgepath $k/(k+l)\angleb{p/q}+l/(k+l)\circleb{p/q}$ is $-l/k$.
Since the curve system $k\anglebb{p/q}+l\circlebb{p/q}$ has $abc$-coordinates $k(1,q-1,p)+l(0,q,p)$,
we have $u=b/(a+b)=1-k/\{(k+l)q\}$.
Then the contribution is calculated by $1-1/\{q(1-u)\}$.

For every subsurface,
the last curve system of the sequence of curve systems corresponding to the subsurface has the common $a$ and $b$ coordinates, say $a_0$ and $b_0$, respectively.
On every hemisphere of $\partial B_i$ divided by the axis,
$(2a_0+b_0)$ subarcs of $\partial S_i$ exist.
$(2a_0+b_0)$ subarcs on each of $2N$ hemisphere are glued each other first, and then $(2a_0+2b_0)$ disks intersecting the axis are connected next.
Thus, the effect of the gluing on the Euler characteristic is $-(2a_0+b_0)\NumTangles+(2a_0+2b_0)$.
Since $\sharp s=a_0$ and $u_0=b_0/(a_0+b_0)$, the ratio of this effect to the number of sheets
is $\{-(2a_0+b_0)\NumTangles+2a_0+2b_0\}/a_0=-1/(1-u_0)\cdot(\NumTangles-2)-\NumTangles$.

  \begin{figure}[htb]
   \begin{center}
    \begin{picture}(242,103)
     \put(0,0){\scalebox{0.5}{\includegraphics{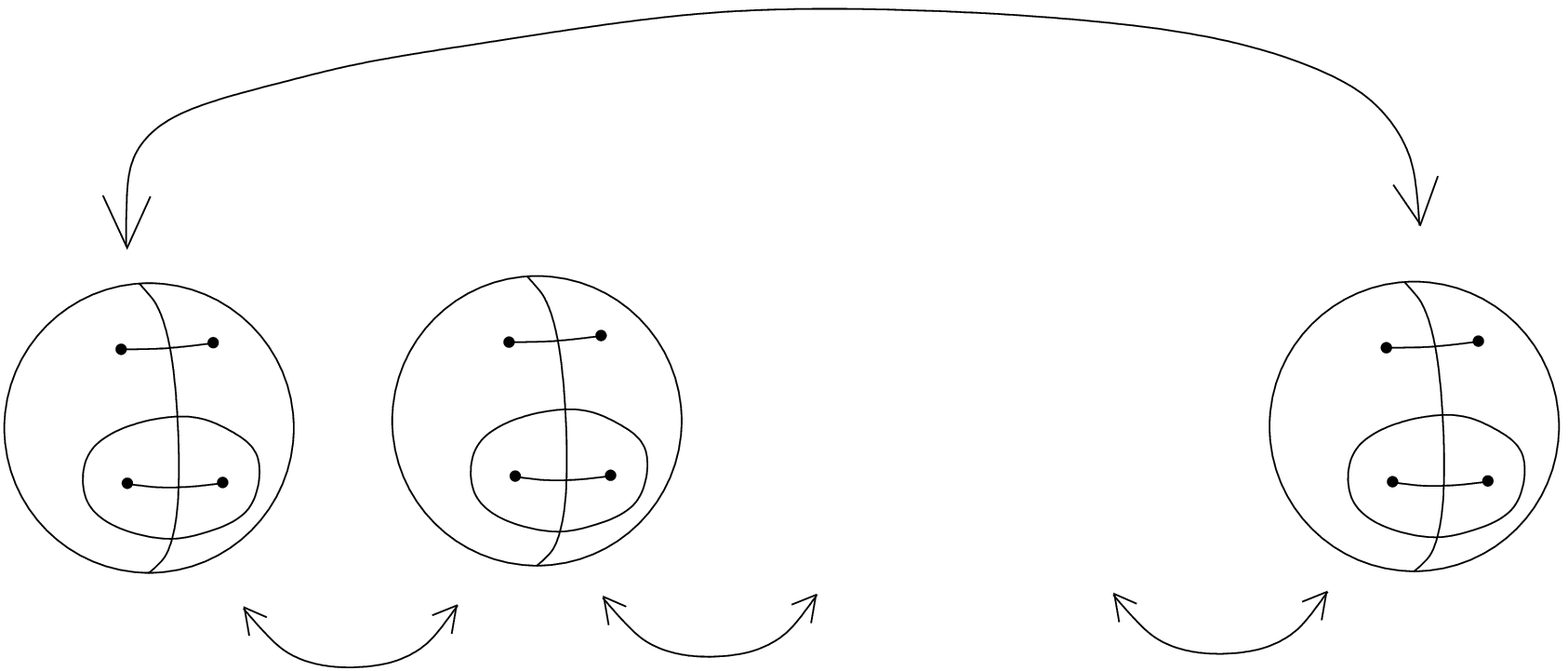}}} 
     \put(18,32){$a_0$}
     \put(33,32){$a_0$}
     \put(18,53){$a_0$}
     \put(31,53){$a_0$}
     \put(7,28){$b_0$}
     \put(133,32){\scalebox{3.0}{$\cdot\cdot\cdot$}}
    \end{picture}
   \end{center}
   \caption{Gluing $\NumTangles$ subsurfaces at subarcs on $\partial B_i$'s}
   \label{Fig:Gluing}
  \end{figure}

Hence,
\begin{eqnarray*}
\chi&=&
\sum_{i=1}^{\NumTangles} (2\cdot \sharp s)
+\sum_{\Edgepath_{i}\in\ConstantEdgepaths} 
 \left(
     (\frac{1}{q_i(1-u)}-1) \cdot \sharp s 
 \right)
\\
&&
-\sum_{\Edgepath_{i}\in\NonconstantEdgepaths} 
 \left( |\Edgepath_{i} | \cdot \sharp s \right)
\\
&&
+\{-\frac{1}{1-u}(\NumTangles-2)-\NumTangles\} \cdot\sharp s
,
\end{eqnarray*}
\begin{eqnarray}
\frac{-\chi}{\sharp s}
&=&
 \sum_{i=1}^{\NumTangles}
  \left(
   \left\{
    \begin{array}{l}
     0 \\
     ~~~~~(\textrm{  if $\Edgepath_i$ is constant }) \\
     |\Edgepath_{i} |\\
     ~~~~~(\textrm{ otherwise }) \\
    \end{array}
   \right.
  \right)
\label{Eq:Formula:EulerCharTypeI}
\\
&& 
+\NumTangles_{\mathrm{const}}-\NumTangles
+\left(
\NumTangles-2-\sum_{\Edgepath_i \in \ConstantEdgepaths}\frac{1}{q_i}
\right)
\frac{1}{1-u}
,
\nonumber
\end{eqnarray}
where the edgepath system $\EdgepathSystem$ is divided into
the set $\NonconstantEdgepaths$ of non-constant edgepaths
and the set $\ConstantEdgepaths$ of constant edgepaths,
and $\NumTangles_{\mathrm{const}}$ denotes the number of the constant edgepaths.



\subsection*{Number of sheets}

When we construct a surface from an edgepath system,
the number of sheets of a surface denoted by $\sharp s$ is determined as follows.

Assume first that the last edge of an edgepath of an edgepath system is a partial edge of length $k/m$ where the fraction $k/m$ is irreducible.
Since the number of saddles $k/m \cdot \sharp s$ must be an integer,
$\sharp s$ is a multiple of $m$.
Assume next that an edgepath of an edgepath system is a constant edgepath 
$k/m \angleb{p/q}$\,--\,$(1-k/m) \circleb{p/q}$.
Since the number of caps $(m-k)/k \cdot \sharp s$ must be an integer,
$\sharp s$ is a multiple of $k$.
Thus, $\sharp s$ is determined as the least common multiple of these integers.


\subsection*{A remark for proofs}

Here, we give an elementary fact for the subsequent sections.
\begin{remark}\label{Rem:Condition:ForLinkToBeKnot}
Let $L$ be a Montesinos link $L(K_1,K_2,\ldots,K_\NumTangles)$
where $\NumTangles\ge 3$ is the number of tangles and 
each $K_i$ is a non-integral rational number.
For the link $L$ to be a knot,
fractions $K_1,K_2,\ldots,K_\NumTangles$ must satisfy either of:
\begin{itemize}
\item
Exactly one of the fractions has even denominator.
\item
All denominators are odd and the number of odd numerators is odd.
\end{itemize}
Moreover, $\sum K_i \neq 0$ holds under the condition.
This means that a type I edgepath system with all its edgepaths being constant 
does not exist for a Montesinos knot.
\end{remark}


%
%

\section{A bound on the denominator}
\label{Sec:ABoundOnTheDenominator}

The purpose of this section is to prove Theorem \ref{Thm:Denom:UpperBound:Main}
about an upper bound of the denominator of a boundary slope.
Proving Lemma \ref{Lem:LowerBound:-XoverS} which claims a lower bound of ${-\chi}/{\sharp s}$,
immediately gives the theorem.
We also show the best possibility and some corollaries.


\subsection{A lower bound of $-\chi/\sharp s$}
\label{SubSec:ALowerBoundOf-XoverS}

This subsection is devoted to giving a proof of Lemma \ref{Lem:LowerBound:-XoverS}, which claims a lower bound of $-\chi/\sharp s$.

\begin{lemma}
\label{Lem:LowerBound:-XoverS}
Let $\Slope=\Numer/\Denom$ be a finite boundary slope.
The corresponding essential surface $F$
has Euler characteristic $\chi$,
the number of sheets $\sharp s$,
and the number of boundary components $\sharp b$.
Then, except for some boundary slopes,
these values satisfy
\[
{-\chi}/{\sharp s}\ge 1.
\]

The exceptions occur from $(-2,3,t)$-pretzel knots for odd $t\ge 3$ or their mirror images.
One of the boundary slopes for the knot only satisfies 
\[
{-\chi}/{\sharp s}\ge (\Denom-1)/\Denom
,
\]
though a slightly stronger condition $\sharp b\ge 2$ holds at the same time.
\end{lemma}

\subsubsection{Type II surfaces and Type III surfaces}

A lower bound of $-\chi/\sharp s$ is easily obtained
for type II surfaces and type III surfaces.

\begin{lemma}
$-\chi/\sharp s \ge 1$
holds for any type II surface and any type III surface 
corresponding to any Montesinos knot with $\NumTangles \ge 3$ tangles.
\end{lemma}

\begin{proof}
For a type II surface,
since every edgepath has at least $1$ complete edge in $u>0$,
we have
$
-\chi/\sharp s
=
 (\sum_{i=1}^{\NumTangles}|\Edgepath_{i,>0} |)
 +|\EdgepathSystem(+0) |-2
\ge 1
$.
For a type III surface,
$-\chi/\sharp s
=
\left( 
 \sum_{i=1}^{\NumTangles}|\Edgepath_{i,>0} |
\right)
\ge 3
$.
By the simplification, we can ignore the effect by the augmentation and type III surfaces with partial $\infty$-edges.
%
\end{proof}

\subsubsection{Type I surfaces}


For type I surfaces, to verify the bound is not so easy as type II and type III surfaces.
Though, for a major part of the type I surfaces,
the bound is shown by ``denominator sequences'' only.

Fix a Montesinos knot $K$ and 
a basic edgepath system $\BasicEdgepathSystem$.
Though type I edgepath systems must correspond to the solution $u$ of the equation (\ref{Eq:EquationForEdgepathSystem}),
we can formally calculate $-\chi/\sharp s$ by the formula (\ref{Eq:Formula:EulerCharTypeI}) for arbitrary $0<u<1$.
Thus, we have a function $X_{\BasicEdgepathSystem}(u):(0,1)\rightarrow \mathbb{R}$.
The function depends on the basic edgepath system $\BasicEdgepathSystem$. 
Though, by examining the formula (\ref{Eq:Formula:EulerCharTypeI})
together with the formula (\ref{Eq:Formula:LengthOfPartialEdge}) about lengths of partial edges,
we can confirm that the function $X_{\BasicEdgepathSystem}$ does not depend on the numerators of vertices which edgepaths pass through or reach.
Namely, for a basic edgepath $\BasicEdgepath_i=\angleb{p_{i,j}/q_{i,j}}$\,--\,$\angleb{p_{i,j-1}/q_{i,j-1}}$\,--\,$\cdots$\,--\,$\angleb{p_{i,2}/q_{i,2}}$\,--\,$\angleb{p_{i,1}/q_{i,1}}$,
$X_{\BasicEdgepathSystem}$ depends on only a $N$-tuple of {\em sequences of denominators} of the form $q_{i,j}(=1)$\,--\,$q_{i,j-1}$\,--\,$\cdots$\,--\,$q_{i,2}$\,--\,$q_{i,1}$.

We introduce a preorder of basic edgepaths and basic edgepath systems.
For two basic edgepaths $\BasicEdgepath_a$ and $\BasicEdgepath_b$,
we say that $\BasicEdgepath_a \le \BasicEdgepath_b$
if $q_{a,k}\le q_{b,k}$ for all $k=1,2,\ldots,\min(j_a,j_b)$
where $j_a$ and $j_b$ mean the lengths of their denominator sequences. 
For two basic edgepath systems $\BasicEdgepathSystem_a$ and $\BasicEdgepathSystem_b$,
we define a preorder by $\BasicEdgepathSystem_a \le \BasicEdgepathSystem_b$ if $\BasicEdgepath_{a,i} \le \BasicEdgepath_{b,i}$ is satisfied for all indices $i$.
It is easy to confirm that
if $\BasicEdgepathSystem_a\le \BasicEdgepathSystem_b$,
then $X_{\BasicEdgepathSystem_a}(u)\le X_{\BasicEdgepathSystem_b}(u)$ holds for any $u \in (0,1)$.

By elementary calculations,
for a basic edgepath system $\BasicEdgepathSystem$
whose set of denominator sequences is one of
\begin{eqnarray}
\begin{minipage}{10cm}
\begin{itemize}
\item \{$1$\,--\,$2$, $1$\,--\,$3$, $1$\,--\,$3$, $1$\,--\,$3$, $\ldots$\} ($\NumTangles \ge 4$), 
\item \{$1$\,--\,$2$, $1$\,--\,$7$, $1$\,--\,$3$\,--\,$7$\},
\item \{$1$\,--\,$2$, $1$\,--\,$2$\,--\,$7$, $1$\,--\,$2$\,--\,$7$\},
\item \{$1$\,--\,$3$, $1$\,--\,$4$, $1$\,--\,$7$\},
\item \{$1$\,--\,$3$, $1$\,--\,$4$, $1$\,--\,$2$\,--\,$5$\},
\item \{$1$\,--\,$3$, $1$\,--\,$5$, $1$\,--\,$5$\},
\item \{$1$\,--\,$4$, $1$\,--\,$4$, $1$\,--\,$4$\},
\end{itemize}
\end{minipage}
\label{Items:DenomSequences}
\end{eqnarray}
the inequality $X_{\BasicEdgepathSystem}(u)\ge 1$ holds for arbitrary $u$.
For instance,
if a basic edgepath system $\BasicEdgepathSystem$ has the set of denominator sequences \{$1$\,--\,$4$, $1$\,--\,$4$, $1$\,--\,$4$\},
$X_{\BasicEdgepathSystem}(u)$ is $1$ for $0<u\le 3/4$ and is $1/\{4(1-u)\}$ for $3/4\le u<1$.
%
Necessarily, another edgepath system greater than such an edgepath system also satisfies $X_{\BasicEdgepathSystem}(u)\ge 1$.
Note that \{$1$\,--\,$2$, $1$\,--\,$3$, $1$\,--\,$3$, $1$\,--\,$3$, \ldots \} ($\NumTangles\ge4$) is the denominator sequence of one of the smallest basic edgepath systems for fixed $\NumTangles$.

Hence, in a sense, a major part of the edgepath systems corresponding to candidate surfaces satisfies $-\chi/\sharp s \ge 1$.
We have:
\begin{lemma}
\label{Lem:LowerBound:-XoverS:B}
Assume that a type I edgepath system $\EdgepathSystem$ 
is included in the extended basic edgepath system of a basic edgepath system $\BasicEdgepathSystem$.
If the basic edgepath system $\BasicEdgepathSystem$
is equal to or greater than one of 
the basic edgepath systems listed in (\ref{Items:DenomSequences}).
Then $\EdgepathSystem$ satisfies $-\chi/\sharp s\ge 1$.
Especially,
any type I surface of a Montesinos knot with $\NumTangles\ge4$ tangles 
always satisfies the inequality.
\end{lemma}


\subsubsection{Remaining cases}

We only have to check for the rest of the edgepath systems.
For the remaining basic edgepath systems $\BasicEdgepathSystem$,
we concretely solve the equation $\sum \widetilde{\BasicEdgepath_i}(u)=0$,
enumerate all the candidate edgepath systems,
and calculate $-\chi/\sharp s$ one by one.
%
Remaining cases are described by denominator sequences as follows.
\begin{eqnarray}
\begin{minipage}{10cm}
\begin{itemize}
 \item $\BasicEdgepath_1:1$--$2$, $\BasicEdgepath_2:\cdots$\,--$3$, $\BasicEdgepath_3:$ arbitrary,
 \item $\BasicEdgepath_1:1$--$2$, $\BasicEdgepath_2:\cdots$\,--$5$, $\BasicEdgepath_3:$ arbitrary,
 \item $\BasicEdgepath_1:1$--$2$, $\BasicEdgepath_2:1$--$t_1$, $\BasicEdgepath_3:1$--$t_2$ ($t_1,t_2$ odd and $\ge7$),
 \item $\BasicEdgepath_1:1$--$2$, $\BasicEdgepath_2:1$--$t_1$, $\BasicEdgepath_3:1$--$2$--$t_2$ ($t_1,t_2$ odd and $\ge7$),
 \item $\BasicEdgepath_1:1$--$3$ or $1$--$2$--$3$, $\BasicEdgepath_2:1$--$3$ or $1$--$2$--$3$, $\BasicEdgepath_3:$ arbitrary,
 \item $\BasicEdgepath_1:1$--$3$ or $1$--$2$--$3$, $\BasicEdgepath_2:1$--$4$ or $\cdots$\,--$3$--$4$, $\BasicEdgepath_3:1$--$5$.
\end{itemize}
\end{minipage}
\label{Items:RemainingCases}
\end{eqnarray}
Basically, $\BasicEdgepath_1$, $\BasicEdgepath_2$ and $\BasicEdgepath_3$ are arranged
so that the denominators of their starting points $K_i$ are in ascending order.

\begin{lemma}
\label{Lem:LowerBound:RemainingCase}
For the rest of type I surfaces,
that is, the type I surfaces to which Lemma \ref{Lem:LowerBound:-XoverS:B} cannot be applied,
nevertheless Lemma \ref{Lem:LowerBound:-XoverS} holds.
Namely,
except for some boundary slopes,
\[
-\chi/\sharp s\ge 1
\]
holds.
The exceptions occur from $(-2,3,t)$-pretzel knots for odd $t\ge 3$ or their mirror images.
One of the boundary slopes for the knots only satisfies
\[
-\chi/\sharp s\ge (\Denom-1)/\Denom
\]
though a slightly stronger condition $\sharp b\ge 2$ holds at the same time.
\end{lemma}

\begin{proof}

First,
we enumerate all choices of a pair ($\BasicEdgepath_1$, $\BasicEdgepath_2$) of basic edgepath systems
whose denominator sequences are included in the list (\ref{Items:RemainingCases}).
Without loss of generality,
we can normalize edgepath systems 
by assuming that the tangles $K_1$ and $K_2$ of the Montesinos knot $K$ satisfy $0<K_1,K_2<1$
and that the last edge of the edgepath $\BasicEdgepath_1$ is decreasing.
There are 27 possible pairs of two edgepaths as listed in the rest of this subsection.
For each choice of $\BasicEdgepath_1$ and $\BasicEdgepath_2$,
we think about the sum $\BasicEdgepath_1+\BasicEdgepath_2$,
which is a function defined by $\BasicEdgepath_1(u)+\BasicEdgepath_2(u)$.
Then, we seek all edges of $\Diagram$ intersecting the sum 
and all vertices on the sum.
For each such an edge or vertex $\Edge$,
we take its mirror image $\Edge_3=-\Edge$ with respect to the $u$-axis,
and then make a basic edgepath $\BasicEdgepath_3$ including $\Edge_3$.
The triple $(\BasicEdgepath_1, \BasicEdgepath_2, \BasicEdgepath_3)$ is ignored if it does not match with the condition in the list (\ref{Items:RemainingCases}).
$\BasicEdgepath_1$, $\BasicEdgepath_2$ and $\BasicEdgepath_3$
have a solution of the equation (\ref{Eq:EquationForEdgepathSystem})
at the $u$-coordinate of the intersection point.
From the solution, we cut the basic edgepaths ($\BasicEdgepath_1$, $\BasicEdgepath_2$, $\BasicEdgepath_3$)
and obtain constant or non-constant edgepaths ($\Edgepath_1$, $\Edgepath_2$, $\Edgepath_3$).
Note that there are many choices of $\BasicEdgepath_3$ which share $\Edge_3$ as the common tail part but have different parts.
Hence, we implicitly discuss many choices of edgepath $\Edgepath_3$ at the same time, though the edgepath with minimum $-\chi/\sharp s$ is important.

In the detailed argument,
note that the denominator $\Denom$ of the slope is the same as the denominator of the twist.
Besides, 
for a constant edgepath $\Edgepath_i$,
a fraction $k_i/m_i$ denotes a particular ratio
which appears in the description
$((k_i/m_i)\angleb{p_i/q_i}+(1-k_i/m_i)\circleb{p_i/q_i})$
of the unique point of the constant edgepath,
especially when we calculate $\sharp s$.
The following fact is used often in the argument.

\begin{remark}\label{Rem:Condition:ForLinkToBeKnot:2}
The starting points of an edgepath system
correspond to the tangles of the Montesinos link.
Therefore, the starting points must satisfy a condition in Remark \ref{Rem:Condition:ForLinkToBeKnot},
for the corresponding Montesinos link to be a knot.
If an edgepath system does not satisfy the condition,
we must add at least one edge to the beginning of some edgepath of an edgepath system.
\end{remark}

Here, we briefly show the calculation of $-\chi/\sharp s$ for the 27 cases.

\begin{enumerate}


\item {$\BasicEdgepath_1=\angleb{0}$--$\angleb{1/2}$, $\BasicEdgepath_2=\angleb{0}$--$\angleb{1/3}$}.
\label{Item:0-1/2,0-1/3}
In Figure \ref{Fig:Edgepaths},
the left and the middle pictures illustrate extended basic edgepath systems $\widetilde{\BasicEdgepath_1}$ and $\widetilde{\BasicEdgepath_2}$,
while the right figure shows $\widetilde{\BasicEdgepath_1}+\widetilde{\BasicEdgepath_2}$ and possible choices of the edge $\Edge$.

\begin{figure}[hbt]
\begin{picture}(256,144)
 \scalebox{0.8}{
   \put(0,0){
    \put(0,0){\scalebox{1.0}{\includegraphics{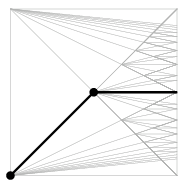}}}
    \put(90,2){u}
    \put(81,-5){1}
    \put(-5,0){0}
    \put(3,-5){0}
    \put(-3,81){1}
    \put(3,91){v}
   }
   \put(110,0){
    \put(0,0){\scalebox{1.0}{\includegraphics{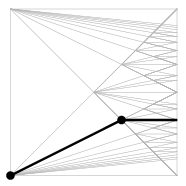}}}
    \put(90,2){u}
    \put(81,-5){1}
    \put(-5,0){0}
    \put(3,-5){0}
    \put(-3,81){1}
    \put(3,91){v}
   }
   \put(220,0){
    \put(0,0){\scalebox{1.0}{\includegraphics{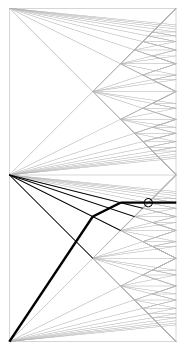}}}
    \put(90,2){u}
    \put(81,-5){1}
    \put(-5,0){0}
    \put(3,-5){0}
    \put(-3,81){1}
    \put(-3,161){2}
    \put(3,169){v}
   }
 }
\end{picture}
\caption{$\widetilde{\BasicEdgepath_1}$, $\widetilde{\BasicEdgepath_2}$, $\widetilde{\BasicEdgepath_1}$+$\widetilde{\BasicEdgepath_2}$, $\Edge$}
\label{Fig:Edgepaths}
\end{figure}

\begin{enumerate}

 \item $\Edge_3=\angleb{-1}$--$\angleb{-1/2}$. 
  \label{Item:0-1/2,0-1/3,A}
   The equation (\ref{Eq:EquationForEdgepathSystem}) gives $3/2\, u=1-u$,
   and we have $u=2/5 ~~(0<2/5<1/2)$.
   From formulae (\ref{Eq:Formula:LengthOfPartialEdge}) and 
   (\ref{Eq:Formula:EulerCharTypeI}),
   $|\Edgepath_1|=\{1+2(2/5 -1)\}/\{(2-1)(2/5 -1)\}
    =1/3$,
   $|\Edgepath_2|=\{1+3(2/5 -1)\}/\{(3-1)(2/5 -1)\}
    =2/3$,
   $|\Edgepath_3|=|\Edgepath_1|
    =1/3$,
   $-\chi/\sharp s =1/3+2/3+1/3+(0-3)+(3-2)\cdot 5/3=4/3 -3+5/3 
    =0$.
  Though, as Remark \ref{Rem:Condition:ForLinkToBeKnot:2}, we must extend the edgepath $\Edgepath_3$ by adding at least one edge. Thus, $-\chi/\sharp s \ge 1$.

 \item $\Edge_3=\angleb{-1}$--$\angleb{-2/3}$. 
  \label{Item:0-1/2,0-1/3,B}
   $u=1/2$,
   $|\Edgepath_1|=0$,
   $|\Edgepath_2|
    =1/2$,
   $|\Edgepath_3|
    =1/2 $,
   $-\chi/\sharp s 
    =0$.
  This example is obtained for the torus knot $K(-1/2,1/3,1/3)$.
  $|\Twist| 
    =2$.
  $\Denom=1$.
  Thus, $-\chi/\sharp s =0=(\Denom-1)/\Denom$.
  Moreover, $\sharp s=\mathrm{lcm}(1,2,2)=2$ and $\Denom=1$ give $\sharp b=2$.

 \item $\Edge_3=\angleb{-1}$--$\angleb{-3/4}$. 
  \label{Item:0-1/2,0-1/3,C}
   $u=3/5  ~~(1/2 <3/5 <2/3 )$,
   $\Edgepath_1$ is constant,
   $|\Edgepath_2|
     =1/4 $,
   $|\Edgepath_3|
     =1/2 $,
   $-\chi/\sharp s 
     =0$.
  By Remark \ref{Rem:Condition:ForLinkToBeKnot:2}, $-\chi/\sharp s \ge 1$.

 \item $\Edge_3=\angleb{-1}$--$\angleb{-4/5}$. 
  \label{Item:0-1/2,0-1/3,D}
   $u=2/3 $,
   $\Edgepath_1$ is constant,
   $|\Edgepath_2|=0$,
   $|\Edgepath_3|
     =1/2 $,
   $-\chi/\sharp s 
      =0$.
  This example is obtained for the torus knot $K(-1/2 ,1/3 ,1/5)$.
    $|\Twist| 
       =1$.
    $\Denom=1$.
    $-\chi/\sharp s =0=\Denom/(\Denom-1)$.
  Moreover, $k_1/m_1 =2/3 $, $\sharp s=\mathrm{lcm}(2,1,2)=2$ and $\Denom=1$ give $\sharp b=2$.

 \item $\Edge_3=\angleb{-5/6}$. 
  \label{Item:0-1/2,0-1/3,E}
   $u=5/6 $,
   $\Edgepath_1$ and
   $\Edgepath_2$ is constant,
   $|\Edgepath_3|=0$,
   $-\chi/\sharp s 
     =0$.
  By Remark \ref{Rem:Condition:ForLinkToBeKnot:2}, $-\chi/\sharp s \ge 1$.
\end{enumerate}


\item {$\BasicEdgepath_1=\angleb{0}$--$\angleb{1/2}$, $\BasicEdgepath_2=\angleb{1}$--$\angleb{1/2}$--$\angleb{1/3}$}.
\label{Item:0-1/2,1-1/2-1/3}
\begin{enumerate}
 \item $\Edge_3=\angleb{-1}$--$\angleb{-(t-1)/t}$ ($t\ge 5$). 
  \label{Item:0-1/2,1-1/2-1/3,B}
   $u=(t-1)/\{2(t-2)\} ~~(1/2 <u\le 2/3 )$,
   $\Edgepath_1$ is constant,
   $|\Edgepath_2|
     =(t-5)/(t-3)$,
   $|\Edgepath_3|
     =(t-4)/(t-3)$,
   $-\chi/\sharp s
     =1-2/(t-3)$.
  \\
    If $t=6$, $-\chi/\sharp s = 1/3 $. 
     By Remark \ref{Rem:Condition:ForLinkToBeKnot:2}, $-\chi/\sharp s = 1/3 +1\ge 1$. 
     If $t$ is even, $-\chi/\sharp s \ge 1$ similarly.
    If $t=5$, $-\chi/\sharp s =0$.
     This example is obtained for the torus knot $K(-1/2 ,1/3 ,1/5 )$.
     For odd $t\ge 7$, 
      $|\Twist|
       =2/(t-3)$.
     $\Denom=(t-3)/2$.
     $(\Denom-1)/\Denom=1-2/(t-3)=-\chi/\sharp s $.
     Moreover, for the constant edgepath $\Edgepath_1$,
     $k_1/m_1 =q_1(1-u)=(t-3)/(t-2)$, $\sharp s=\mathrm{lcm}(t-3,(t-3)/2,t-3)=t-3$, and $\Denom=(t-3)/2$ give $\sharp b=2$.

 \item $\Edge_3=\angleb{-5/6}$. 
   The edgepath system is same as in the item \ref{Item:0-1/2,0-1/3,E}.
\end{enumerate}


\item {$\BasicEdgepath_1=\angleb{0}$--$\angleb{1/2}$, $\BasicEdgepath_2=\angleb{1}$--$\angleb{2/3}$}.
\label{Item:0-1/2,1-2/3}
\begin{enumerate}

 \item $\Edge_3=\angleb{-1}$--$\angleb{-4/3}$. 
  \label{Item:0-1/2,1-2/3,A}
   $u$ is non-isolated ($0\le u\le 1/2 $),
   $|\Edgepath_1|
     =(1-2u)/(1-u)$,
   $|\Edgepath_2|
     =(2-3u)/(2-2u)$,
   $|\Edgepath_3|
     =(2-3u)/(2-2u)$,
   $-\chi/\sharp s 
     =2-1/(1-u)$.
  At $u=0$, $-\chi/\sharp s =1$.
  At $u=1/2 $, $-\chi/\sharp s =0$. 
     $0\le -\chi/\sharp s \le 1$.

  This edgepath system is obtained for the torus knot $K(-1/2 ,1/3 ,1/3 )$.

    $|\Twist|
       =2$.
    $\Denom=1$.
  $(\Denom-1)/\Denom=0\le -\chi/\sharp s $.

    If $0<u<1/2 $, $|\Edgepath_1|\notin\mathbb{Z}$ gives $\sharp s\ge 2$. Since $\Denom=1$, we have $\sharp b\ge 2$.
    In the case $u=0$, the surface can be regarded as Type II.
    If $u=1/2 $, 
    the argument reduces to the item \ref{Item:0-1/2,0-1/3,B} corresponding to the torus knot $K(-1/2 ,1/3 ,1/3 )$.

 \item $\Edge_3=\angleb{-1}$--$\angleb{-5/4}$.
  \label{Item:0-1/2,1-2/3,B}
  The calculation is same as in the item \ref{Item:0-1/2,0-1/3,C}, though the edgepath systems themselves do not coincide with each other.

 \item $\Edge_3=\angleb{-1}$--$\angleb{-6/5}$. 
  \label{Item:0-1/2,1-2/3,C}
  The calculation is same as in the item \ref{Item:0-1/2,0-1/3,D}.

 \item $\Edge_3=\angleb{-7/6}$. 
  The calculation is same as in the item \ref{Item:0-1/2,0-1/3,E}.

\end{enumerate}


\item {$\BasicEdgepath_1=\angleb{0}$--$\angleb{1/2}$, $\BasicEdgepath_2=\angleb{0}$--$\angleb{1/2}$--$\angleb{2/3}$}.
\label{Item:0-1/2,0-1/2-2/3}
\begin{enumerate}

 \item $\Edge_3=\angleb{-1}$--$\angleb{-(t-1)/t}$ ($t\ge 2$). 
  \label{Item:0-1/2,0-1/2-2/3,A}
   $u=(t-1)/(2t-1) ~~(1/3 \le u<1/2 )$,
   $|\Edgepath_1|
     =1/t $,
   $|\Edgepath_2|
     =1+1/t $,
   $|\Edgepath_3|
     =(t-1)/t$,
   $-\chi/\sharp s 
     =1$.

 \item $\Edge_3=\angleb{-1}$--$\angleb{-(t+1)/t}$ ($t\ge 5$).
  \label{Item:0-1/2,0-1/2-2/3,B}
  The calculation is same as in the item \ref{Item:0-1/2,1-1/2-1/3,B}.

 \item $\Edge_3=\angleb{-7/6}$.
  The calculation is same as in the item \ref{Item:0-1/2,0-1/3,E}.

\end{enumerate}


\item {$\BasicEdgepath_1=\angleb{0}$--$\angleb{1/2}$, $\BasicEdgepath_2=\angleb{0}$--$\angleb{1/5}$}.
\label{Item:0-1/2,0-1/5}
\begin{enumerate}
 \item $\Edge_3=\angleb{-1}$--$\angleb{-1/2}$. 
  \label{Item:0-1/2,0-1/5,A}
   $u=4/9  ~~(0<4/9 <1/2 )$,
   $|\Edgepath_1|
     =1/5$,
   $|\Edgepath_2|
     =4/5$,
   $|\Edgepath_3|
     =1/5 $,
   $-\chi/\sharp s 
     =0$.
   By Remark \ref{Rem:Condition:ForLinkToBeKnot:2}, $-\chi/\sharp s \ge 1$.

 \item $\Edge_3=\angleb{-2/3}$. 
  \label{Item:0-1/2,0-1/5,B}
   $u=2/3$,
   $\Edgepath_1$ is constant,
   $|\Edgepath_2|
     =1/2 $,
   $|\Edgepath_3|=0$,
   $-\chi/\sharp s 
     =0$.
  For the denominators of tangles to be in ascending order, we must add at least an edge to the edgepath $\Edgepath_3$. Hence, $-\chi/\sharp s \ge 1$.

 \item $\Edge_3=\angleb{-2/3}$--$\angleb{-5/7}$. 
  \label{Item:0-1/2,0-1/5,C}
   $u$ is non-isolated ($2/3 \le u\le 4/5 $),
   $\Edgepath_1$ is constant,
   $|\Edgepath_2|
     =(4-5u)/(4-4u)$, 
   $|\Edgepath_3|
     =(6-7u)/(4-4u)$,
   $-\chi/\sharp s 
     =1$.

 \item $\Edge_3=\angleb{-7/10}$. 
  \label{Item:0-1/2,0-1/5,D}
   $u=9/10 $,
   $\Edgepath_1$ and 
   $\Edgepath_2$ is constant,
   $|\Edgepath_3|=0$,
   $-\chi/\sharp s 
     =2$.

\end{enumerate}


\item {$\BasicEdgepath_1=\angleb{0}$--$\angleb{1/2}$, $\BasicEdgepath_2=\angleb{1}$--$\angleb{1/2}$--$\angleb{1/3}$--$\angleb{1/4}$--$\angleb{1/5}$}.
\label{Item:0-1-2,1-1/2-1/3-1/4-1/5}
\begin{enumerate}

 \item $\Edge_3=\angleb{-1}$--$\angleb{-(t-1)/t}$ ($t\ge 5$). 
  \label{Item:0-1-2,1-1/2-1/3-1/4-1/5,B}
   $u=(t-1)/\{2(t-2)\} ~(1/2 <u\le 2/3 )$,
   $\Edgepath_1$ is constant,
   $|\Edgepath_2|\ge 2$,
   $|\Edgepath_3|\ge 0$,
   $-\chi/\sharp s \ge 
     1$.

 \item $\Edge_3=\angleb{-3/4}$. 
  \label{Item:0-1-2,1-1/2-1/3-1/4-1/5,A}
   $u=3/4$,
   $\Edgepath_1$ is constant,
   $|\Edgepath_2|=1$,
   $|\Edgepath_3|=0$,
   $-\chi/\sharp s 
     =1$.

 \item $\Edge_3=\angleb{-2/3}$--$\angleb{-5/7}$. 
  \label{Item:0-1-2,1-1/2-1/3-1/4-1/5,C}
   $u=4/5$,
   $\Edgepath_1$ is constant,
   $|\Edgepath_2|=0$,
   $|\Edgepath_3|
     =1/2$,
   $-\chi/\sharp s 
     =1$.

 \item $\Edge_3=\angleb{-7/10}$.
  The edgepath system is same as in the item \ref{Item:0-1/2,0-1/5,D}.
\end{enumerate}


\item {$\BasicEdgepath_1=\angleb{0}$--$\angleb{1/2}$, $\BasicEdgepath_2=\angleb{1}$--$\angleb{1/2}$--$\angleb{2/5}$}.
\label{Item:0-1/2,1-1/2-2/5}
\begin{enumerate}

 \item $\Edge_3=\angleb{-1}$--$\angleb{-(t-1)/t}$ ($t\ge 9$). 
  \label{Item:0-1/2,1-1/2-2/5,A}
   $u=(t-1)/\{2(t-4)\} ~(1/2 <u\le 4/5 )$,
   $\Edgepath_1$ is constant,
   $|\Edgepath_2|
     =(t-9)/(t-7)$,
   $|\Edgepath_3|
     =(t-8)/(t-7)$,
   $-\chi/\sharp s 
     =1$.

 \item $\Edge_3=\angleb{-9/10}$.
  \label{Item:0-1/2,1-1/2-2/5,B}
  The calculation is same as in the item \ref{Item:0-1/2,0-1/5,D}.

\end{enumerate}


\item {$\BasicEdgepath_1=\angleb{0}$--$\angleb{1/2}$, $\BasicEdgepath_2=\angleb{0}$--$\angleb{1/2}$--$\angleb{2/5}$}.
\label{Item:0-1/2,0-1/2-2/5}
\begin{enumerate}

 \item $\Edge_3=\angleb{-1}$--$\angleb{-(t-1)/t}$ ($t\ge 2$) ($0<u<1/2 $). 
  \label{Item:0-1/2,0-1/2-2/5,A}
   $u=(t-1)/(2t-1) ~(1/3 \le u<1/2 )$,
   $|\Edgepath_1|
     =1/t$,
   $|\Edgepath_2|
     =1+1/t $,
   $|\Edgepath_3|
     =(t-1)/t$,
   $-\chi/\sharp s 
     =1$.

 \item $\Edge_3=\angleb{-1}$--$\angleb{-(t-1)/t}$ ($t\ge 9$) ($u>1/2 $).
  \label{Item:0-1/2,0-1/2-2/5,B}
  The edgepath system is same as in the item \ref{Item:0-1/2,1-1/2-2/5,A}.

 \item $\Edge_3=\angleb{-9/10}$.
  The calculation is same as in the item \ref{Item:0-1/2,0-1/5,D}.

\end{enumerate}


\item {$\BasicEdgepath_1=\angleb{0}$--$\angleb{1/2}$, $\BasicEdgepath_2=\angleb{0}$--$\angleb{1/3}$--$\angleb{2/5}$}.
\label{Item:0-1/2,0-1/3-2/5}
\begin{enumerate}

 \item $\Edge_3=\angleb{-1}$--$\angleb{-1/2}$. 
   $u=2/5  ~~(0<2/5 <1/2 )$,
   $|\Edgepath_1|
     =1/3$,
   $|\Edgepath_2|
     =1+2/3$,
   $|\Edgepath_3|
     =1/3 $,
   $-\chi/\sharp s 
     =1$.

 \item $\Edge_3=\angleb{-1}$--$\angleb{-2/3}$. 
   $u=1/2 $,
   $|\Edgepath_1|=0$,
   $|\Edgepath_2|
    =1+1/2$,
   $|\Edgepath_3|=1/2 $,
   $-\chi/\sharp s 
     =1$.

 \item $\Edge_3=\angleb{-1}$--$\angleb{-3/4}$. 
  \label{Item:0-1/2,0-1/3-2/5,C}
   $u=3/5  ~~(1/2 <3/5 <2/3 )$,
   $\Edgepath_1$ is constant,
   $|\Edgepath_2|
     =1+1/4$,
   $|\Edgepath_3|
     =1/2$,
   $-\chi/\sharp s 
     =1$.

 \item $\Edge_3=\angleb{-1}$--$\angleb{-(t-1)/t}$ ($5\le t\le 9$). 
  \label{Item:0-1/2,0-1/3-2/5,D}
   $u=(t-1)/(t+1) ~~(2/3 \le u \le 4/5 )$,
   $\Edgepath_1$ is constant,
   $|\Edgepath_2|
     =(9-t)/4$,
   $|\Edgepath_3|
     =1/2$,
   $-\chi/\sharp s 
     =1$.

 \item $\Edge_3=\angleb{-9/10}$.
  The calculation is same as in the item \ref{Item:0-1/2,0-1/5,D}.

\end{enumerate}


\item {$\BasicEdgepath_1=\angleb{0}$--$\angleb{1/2}$, $\BasicEdgepath_2=\angleb{1}$--$\angleb{1/2}$--$\angleb{3/5}$}.
\begin{enumerate}

 \item $\Edge_3=\angleb{-1}$--$\angleb{-(t+1)/t}$ ($t\ge 9$).
  The calculation is same as in the item \ref{Item:0-1/2,1-1/2-2/5,A}.

 \item $\Edge_3=\angleb{-11/10}$.
  The calculation is same as in the item \ref{Item:0-1/2,0-1/5,D}.

\end{enumerate}


\item {$\BasicEdgepath_1=\angleb{0}$--$\angleb{1/2}$, $\BasicEdgepath_2=\angleb{0}$--$\angleb{1/2}$--$\angleb{3/5}$}.
\begin{enumerate}

 \item $\Edge_3=\angleb{-1}$--$\angleb{-(t-1)/t}$.
  The calculation is same as in the item \ref{Item:0-1/2,0-1/2-2/5,A}.

 \item $\Edge_3=\angleb{-1}$--$\angleb{-(t+1)/t}$ ($t\ge 9$).
  The calculation is same as in the item \ref{Item:0-1/2,1-1/2-2/5,A}.

 \item $\Edge_3=\angleb{-11/10}$.
  The calculation is same as in the item \ref{Item:0-1/2,0-1/5,D}.

\end{enumerate}


\item {$\BasicEdgepath_1=\angleb{0}$--$\angleb{1/2}$, $\BasicEdgepath_2=\angleb{1}$--$\angleb{2/3}$--$\angleb{3/5}$}.
\begin{enumerate}

 \item $\Edge_3=\angleb{-1}$--$\angleb{-4/3}$. 
   $u$ is non-isolated ($0\le u\le 1/2 $),
   $|\Edgepath_1|
     =(1-2u)/(1-u)$,
   $|\Edgepath_3|
     =(2-3u)/(2-2u)$,
   $|\Edgepath_2|
     =1+(2-3u)/(2-2u)$,
   $-\chi/\sharp s 
     =1$.

 \item $\Edge_3=\angleb{-1}$--$\angleb{-5/4}$.
  The calculation is same as in the item \ref{Item:0-1/2,0-1/3-2/5,C}.

 \item $\Edge_3=\angleb{-1}$--$\angleb{-(t+1)/t}$ ($5\le t \le 9$).
  The calculation is same as in the item \ref{Item:0-1/2,0-1/3-2/5,D}.

 \item $\Edge_3=\angleb{-11/10}$.
  The calculation is same as in the item \ref{Item:0-1/2,0-1/5,D}.
\end{enumerate}


\item {$\BasicEdgepath_1=\angleb{0}$--$\angleb{1/2}$, $\BasicEdgepath_2=\angleb{1}$--$\angleb{4/5}$}.
\begin{enumerate}
 \item $\Edge_3=\angleb{-4/3}$.
  The calculation is same as in the item \ref{Item:0-1/2,0-1/5,B}.
 \item $\Edge_3=\angleb{-4/3}$--$\angleb{-9/7}$.
  The calculation is same as in the item \ref{Item:0-1/2,0-1/5,C}.
 \item $\Edge_3=\angleb{-13/10}$.
  The calculation is same as in the item \ref{Item:0-1/2,0-1/5,D}.
\end{enumerate}


\item {$\BasicEdgepath_1=\angleb{0}$--$\angleb{1/2}$, $\BasicEdgepath_2=\angleb{0}$--$\angleb{1/2}$--$\angleb{2/3}$--$\angleb{3/4}$--$\angleb{4/5}$}.
\begin{enumerate}
 \item $\Edge_3=\angleb{-1}$--$\angleb{-(t-1)/t}$ ($t\ge2$). 
   $u=(t-1)/(2t-1) ~(1/3 \le u<1/2 )$,
   $|\Edgepath_1|
     =1/t $,
   $|\Edgepath_2|
     =3+1/t$,
   $|\Edgepath_3|
     =(t-1)/t$,
   $-\chi/\sharp s 
     =3$.

 \item $\Edge_3=\angleb{-1}$--$\angleb{-(t+1)/t}$ ($t\ge5$).
  The calculation is same as in the item \ref{Item:0-1-2,1-1/2-1/3-1/4-1/5,B}.
 \item $\Edge_3=\angleb{-1}$--$\angleb{-5/4}$.
  The calculation is same as in the item \ref{Item:0-1-2,1-1/2-1/3-1/4-1/5,A}.
 \item $\Edge_3=\angleb{-4/3}$--$\angleb{-9/7}$.
  The calculation is same as in the item \ref{Item:0-1-2,1-1/2-1/3-1/4-1/5,C}.
 \item $\Edge_3=\angleb{-13/10}$.
  The calculation is same as in the item \ref{Item:0-1/2,0-1/5,D}.
\end{enumerate}


\item {$\BasicEdgepath_1=\angleb{0}$--$\angleb{1/2}$, $\BasicEdgepath_2=\angleb{0}$--$\angleb{1/t}$ (odd $t\ge 7$)}.
\label{Item:0-1/2,0-1/t}
\begin{enumerate}

  Properly speaking, $t \le 2p-1$ is required in the items \ref{Item:0-1/2,0-1/t,B}, \ref{Item:0-1/2,0-1/t,D} and \ref{Item:0-1/2,0-1/t,C} below.
  Though, in this part,
  we allow the denominator of the second tangle to be greater than that of the third tangle exceptionally. 
  Thus, we avoid repeating essentially the same calculations.

 \item $\Edge_3=\angleb{-1}$--$\angleb{-1/2}$. 
  \label{Item:0-1/2,0-1/t,A}
   $u=(t-1)/(2t-1)$,
   $|\Edgepath_1|
     =1/t$,
   $|\Edgepath_2|
     =(t-1)/t$,
   $|\Edgepath_3|
     =1/t$,
   $-\chi/\sharp s 
     =0$.
  By Remark \ref{Rem:Condition:ForLinkToBeKnot:2}, $-\chi/\sharp s \ge 1$.

 \item $\Edge_3=\angleb{-1/2}$--$\angleb{-p/(2p-1)}$ ($p\ge 4$, $u\le (t-1)/t$). 
  \label{Item:0-1/2,0-1/t,B}
  $p\ge 4$ is derived from $2p-1\ge 7$.
  By the conditions
  $1/2+1/(t-1)\cdot u-\{1/(2p-3)\cdot(u-1/2)+1/2\}=0$ from (\ref{Eq:EquationForEdgepathSystem}),
  $1/2\le u\le (t-1)/t$ and
  $1/2\le u\le (2p-2)/(2p-1)$,
  the solution $u=(t-1)/\{2(t-2p+2)\}$ exists if $t-4p+3\ge 0$ holds.
  $\Edgepath_1$ is constant, 
  $|\Edgepath_2|=(t-4p+4)/(t-4p+5)$,
  $|\Edgepath_3|=(t-4p+3)/(t-4p+5)$,
  $-\chi/\sharp s=1+(2p-6)/(t-4p+5)\ge 1$.

 \item $\Edge_3=\angleb{-1/2}$--$\angleb{-p/(2p-1)}$ ($p\ge 2$, $u\ge (t-1)/t$). 
  \label{Item:0-1/2,0-1/t,D}
  By the conditions
  $1/2+1/t-\{1/(2p-3)\cdot(u-1/2)+1/2\}=0$ from (\ref{Eq:EquationForEdgepathSystem}),
  $u\ge (t-1)/t$ and
  $1/2\le u\le (2p-2)/(2p-1)$,
  we have a contradiction $4p-2\le t\le 4p-4$.
  Thus, no solution exists.

 \item $\Edge_3=\angleb{-p/(2p-1)}$--$\circleb{-p/(2p-1)}$ ($p\ge 4$). 
  \label{Item:0-1/2,0-1/t,C}
  By the conditions
  $1/2+1/(t-1)\cdot u-p/(2p-1)=0$,
  $u\ge (2p-2)/(2p-1)$ and
  $1/2\le u\le (t-1)/t$,
  the solution $u=(t-1)/\{2(2p-1)\}$ exists if $4p-3\le t\le 4p-2$ holds.
  The solution for $t=4p-3$ is treated in \ref{Item:0-1/2,0-1/t,B},
  while $t=4p-2$ is even and unsuitable. 

\end{enumerate}


\item {$\BasicEdgepath_1=\angleb{0}$--$\angleb{1/2}$, $\BasicEdgepath_2=\angleb{0}$--$\angleb{(t-1)/t}$ (odd $t\ge 7$)}.
\begin{enumerate}
 \item $\Edge_3=\angleb{-1/2}$--$\angleb{-p/(2p+1)}$ ($p\ge 2$). 
  For the denominators of the tangles to be in ascending order, $7\le t\le 2p-1$ and $p\ge 4$ are required.
  The calculation is similar to the item \ref{Item:0-1/2,0-1/t,B}.
\end{enumerate}


\item {$\BasicEdgepath_1=\angleb{0}$--$\angleb{1/2}$, $\BasicEdgepath_2=\angleb{0}$--$\angleb{1/2}$--$\angleb{p/(2p\pm 1)}$}.
\begin{enumerate}

 \item $\Edge_3=\angleb{-1}$--$\angleb{-1/t}$. 
  The calculation is same as in the item \ref{Item:0-1/2,0-1/t,B}.
\end{enumerate}


\item {$\BasicEdgepath_1=\angleb{0}$--$\angleb{1/3}$, $\BasicEdgepath_2=\angleb{0}$--$\angleb{1/3}$}.
\label{Item:0-1/3,0-1/3}

\begin{enumerate}

 \item $\Edge_3=\angleb{0}$--$\angleb{-1/2}$.
  \label{Item:0-1/3,0-1/3,A}
   %
   $u$ is non-isolated ($0\le u\le 1/2$), 
   $|\Edgepath_1|
     =(3u-2)/(2u-2)$,
   $|\Edgepath_2|
     =(3u-2)/(2u-2)$,
   $|\Edgepath_3|
     =(2u-1)/(u-1)$,
   $-\chi/\sharp s 
     =2-1/(1-u)\ge 0$.
  Similarly to the item \ref{Item:0-1/2,1-2/3,A},
  $-\chi/\sharp s \ge (\Denom-1)/\Denom$ holds.

 \item $\Edge_3=\angleb{-1/2}$--$\angleb{-2/3}$.
   $u$ is non-isolated ($1/2 \le u \le 2/3 $),
   $|\Edgepath_1|
      =(3u-2)/(2u-2)$,
   $|\Edgepath_2|
      =(3u-2)/(2u-2)$,
   $|\Edgepath_3|
      =(3u-2)/(u-1)$,
   $-\chi/\sharp s 
      =3-1/(1-u)\ge 0$.
  By Remark \ref{Rem:Condition:ForLinkToBeKnot:2}, $-\chi/\sharp s \ge 1$.

\end{enumerate}


\item {$\BasicEdgepath_1=\angleb{0}$--$\angleb{1/3}$, $\BasicEdgepath_2=\angleb{1}$--$\angleb{2/3}$}. 
  No solution exists.


\item {$\BasicEdgepath_1=\angleb{0}$--$\angleb{1/3}$, $\BasicEdgepath_2=\angleb{1}$--$\angleb{1/2}$--$\angleb{1/3}$}.
\begin{enumerate}
 \item $\Edge_3=\angleb{-1}$--$\angleb{-2/3}$. 
   $u$ is non-isolated ~~($0\le u\le1/2 $),
   $|\Edgepath_1|
     =(3u-2)/(2u-2)$,
   $|\Edgepath_2|
     =1+(2u-1)/(u-1)$,
   $|\Edgepath_3|
     =(3u-2)/(2u-2)$,
   $-\chi/\sharp s 
     =3-1/(1-u)\ge 1$.

 \item $\Edge_3=\angleb{-1}$--$\angleb{-2/3}$.
   $u$ is non-isolated ~~($1/2 < u\le2/3 $),
   $|\Edgepath_1|
     =(3u-2)/(2u-2)$,
   $|\Edgepath_2|
     =(3u-2)/(u-1)$,
   $|\Edgepath_3|
     =(3u-2)/(2u-2)$,
   $-\chi/\sharp s 
     =3-1/(1-u)\ge 0$.
   By Remark \ref{Rem:Condition:ForLinkToBeKnot:2}, $-\chi/\sharp s \ge 1$.
\end{enumerate}


\item {$\BasicEdgepath_1=\angleb{0}$--$\angleb{1/3}$, $\BasicEdgepath_2=\angleb{0}$--$\angleb{1/2}$--$\angleb{2/3}$}.
\begin{enumerate}

 \item $\Edge_3=\angleb{-1}$--$\angleb{-1/2}$. 
   $u=2/5  ~~(0<2/5 <1/2 )$,
   $|\Edgepath_1|
     =2/3 $,
   $|\Edgepath_2|
     =1+1/3 $,
   $|\Edgepath_3|
     =1/3 $,
   $-\chi/\sharp s 
     =1$.

 \item $\Edge_3=\angleb{-1}$--$\angleb{-(t-1)/t}$ ($t\ge 3$). 
   $u=(2t-2)/(3t-1) ~(1/2 \le u <2/3 )$,
   $|\Edgepath_1|
     =2/(t+1)$,
   $|\Edgepath_2|
     =4/(t+1)$,
   $|\Edgepath_3|
     =(t-1)/(t+1)$,
   $-\chi/\sharp s 
     =1$.
\end{enumerate}


\item {$\BasicEdgepath_1=\angleb{0}$--$\angleb{1/2}$--$\angleb{2/3}$, $\BasicEdgepath_2=\angleb{0}$--$\angleb{1/2}$--$\angleb{2/3}$}.
\begin{enumerate}
 \item $\Edge_3=\angleb{-1}$--$\angleb{-(t-1)/t}$ ~($t\ge 2)$. 
   $u=(t-1)/(2t-1) ~~(1/3 \le u<1/2 )$,
   $|\Edgepath_1|=1+1/t $,
   $|\Edgepath_2|=1+1/t $,
   $|\Edgepath_3|=1-1/t $,
   $-\chi/\sharp s 
     =2$.

 \item $\Edge_3=\angleb{-1}$--$\angleb{-(t+1)/t}$.
   $u=(t-1)/(2t-3) ~~(1/2 < u \le 2/3 )$,
   $|\Edgepath_1|
     =(t-3)/(t-2)$,
   $|\Edgepath_2|
     =(t-3)/(t-2)$,
   $|\Edgepath_3|
     =(t-3)/(t-2)$,
   $-\chi/\sharp s 
     =2-2/(t-2)$.
  If $t=3$, though $-\chi/\sharp s =0$, by Remark \ref{Rem:Condition:ForLinkToBeKnot:2}, we have $-\chi/\sharp s \ge 1$.
  If $t\ge 4$, $-\chi/\sharp s \ge 1$ holds.
\end{enumerate}


\item {$\BasicEdgepath_1=\angleb{0}$--$\angleb{1/2}$--$\angleb{2/3}$, $\BasicEdgepath_2=\angleb{1}$--$\angleb{1/2}$--$\angleb{1/3}$}.
No solutions exist.


\item {$\BasicEdgepath_1=\angleb{0}$--$\angleb{1/3}$ or $\angleb{1}$--$\angleb{1/2}$--$\angleb{1/3}$, $\BasicEdgepath_2=\angleb{0}$--$\angleb{1/4}$}.
The check is necessary only for the case that the denominators for the edge $\Edge_3$ is $1$--$5$.
Neither $\angleb{0}$--$\angleb{-1/5}$ nor $\angleb{-1}$--$\angleb{-4/5}$ intersects the sum $\BasicEdgepath_1+\BasicEdgepath_2$.
Thus, no solution exists.


\item {$\BasicEdgepath_1=\angleb{0}$--$\angleb{1/3}$ or $\angleb{1}$--$\angleb{1/2}$--$\angleb{1/3}$, $\BasicEdgepath_2=\angleb{1}$--$\angleb{3/4}$}.
\begin{enumerate}
 \item $\Edge_3=\angleb{-1}$--$\angleb{-6/5}$. 
   $u=12/19  ~~(1/2 <12/19 <2/3 )$,
   $|\Edgepath_1|
     =2/7$,
   $|\Edgepath_2|
     =3/7$,
   $|\Edgepath_3|
     =4/7$,
   $-\chi/\sharp s 
     =1$.
\end{enumerate}


\item {$\BasicEdgepath_1=\angleb{0}$--$\angleb{1/3}$ or $\angleb{1}$--$\angleb{1/2}$--$\angleb{1/3}$, $\BasicEdgepath_2=\angleb{1}$--$\angleb{1/2}$--$\angleb{1/3}$--$\angleb{1/4}$}.

\begin{enumerate}
 \item $\Edge_3=\angleb{-1}$--$\angleb{-6/5}$. 
   $u=4/9  ~~(0< 4/9 <1/2 )$,
   $|\Edgepath_1|
     =1+1/5$,
   $|\Edgepath_2|=2+1/5 $,
   $|\Edgepath_3|
     =4/5$,
   $-\chi/\sharp s 
     =3$.

 \item $\Edge_3=\angleb{-1}$--$\angleb{-4/5}$. 
   $u=4/7  ~~(1/2 <4/7 <2/3 )$,
   $|\Edgepath_1|
     =2/3$,
   $|\Edgepath_2|=1+2/3 $,
   $|\Edgepath_3|
     =2/3$,
   $-\chi/\sharp s 
     =7/3$.

\end{enumerate}


\item {$\BasicEdgepath_1=\angleb{0}$--$\angleb{1/3}$ or $\angleb{1}$--$\angleb{1/2}$--$\angleb{1/3}$, $\BasicEdgepath_2=\angleb{0}$--$\angleb{1/2}$--$\angleb{2/3}$--$\angleb{3/4}$}.

\begin{enumerate}
 \item $\Edge_3=\angleb{-1}$--$\angleb{-4/5}$. 
   $u=4/7  ~~(1/2<4/7<2/3 )$,
   $|\Edgepath_1|
     =1/3$,
   $|\Edgepath_2|
     =1+2/3$,
   $|\Edgepath_3|
     =2/3$,
   $-\chi/\sharp s 
     =2$.
   \\ 
\end{enumerate}

\end{enumerate}

Thus, in most of the cases, $-\chi/\sharp s \ge 1$ is satisfied.
The cases in which only $-\chi/\sharp s \ge (\Denom-1)/\Denom$ is satisfied are 
items \ref{Item:0-1/2,0-1/3,B},
\ref{Item:0-1/2,0-1/3,D},
\ref{Item:0-1/2,1-1/2-1/3,B},
\ref{Item:0-1/2,1-2/3,A},
\ref{Item:0-1/2,1-2/3,C},
\ref{Item:0-1/2,0-1/2-2/3,B} and
\ref{Item:0-1/3,0-1/3,A}.
Last three items are reduced to 
\ref{Item:0-1/2,0-1/3,D},
\ref{Item:0-1/2,1-1/2-1/3,B} and
\ref{Item:0-1/2,1-2/3,A}
respectively.
In any of these cases,
the knot is essentially a Montesinos knot $K(-1/2 ,1/3 ,1/t )$ ($t\ge 3$, $t$ is odd).
Eventually, the candidate surfaces with $0\le -\chi/\sharp s < 1$ are
\begin{itemize}
\item[(a)] annuli for the torus knots $K(-1/2,1/3,1/3)$ and $K(-1/2,1/3,1/5)$,
\item[(b)] a surface with $1/2 \le -\chi/\sharp s = 1-2/(t-3) < 1$ corresponding to the knot $K(-1/2,1/3,1/t)$ for odd $t\ge 7$,
\item[(c)] a family of surfaces with $0 < -\chi/\sharp s <1$ which corresponds to the non-isolated solutions for $K(-1/2,1/3,1/3)$.
\end{itemize}
Note that surfaces in the family (c) in the above list are compressible, in fact.

\end{proof}



\subsection{Corollaries and best possibility of Theorem \ref{Thm:Denom:UpperBound:Main}}

\subsubsection{Proofs of corollaries}

Once the theorem is shown,
proof of Corollary \ref{Cor:Denom:UpperBound:Main:ByGenus} is straightforward.
Since the argument does not depend on the orientability,
Corollary \ref{Cor:Denom:UpperBound:Main:ByNonOrientableGenus} is also easily obtained.

\begin{proof}[Proof of Corollary \ref{Cor:Denom:UpperBound:Main:ByGenus}]

First, assume that $\sharp b\ge 2$.
A boundary slope and its corresponding surface satisfy
at least the inequality $\Denom \le -\chi/\sharp b + 1$ 
in Theorem \ref{Thm:Denom:UpperBound:Main}.
With a variable $g=(2-\chi-\sharp b)/2$,
if $g\ge 1$, 
we have 
\begin{eqnarray}
\Denom&\le& \frac{-2+2g+\sharp b}{\sharp b}+1=\frac{-2+2g}{\sharp b}+2\le g+1
.
\end{eqnarray}
If $g=0$, we have $\Denom\le 2-2/{\sharp b}<2$, which means $\Denom=1$.
$\Denom\le g+1$ is satisfied also in this case.

Next, assume that $\sharp b=1$.
Then the inequality $\Denom \le -\chi$ in Theorem \ref{Thm:Denom:UpperBound:Main} is satisfied.
Then,
\begin{eqnarray}
\Denom&\le& -2+2g+\sharp b=2g-1.
\end{eqnarray}

By taking maximum of $g+1$ and $2g-1$,
we have $\Denom\le g+1$ for $g=0,1$ and $\Denom\le 2g-1$ for $g\ge 2$.

In the case of $g=0$, equivalently, the surface is planar, 
the inequality $\Denom \le -\chi$ 
in Theorem \ref{Thm:Denom:UpperBound:Main} cannot be satisfied. 
This means that the cases are exceptional, that is, 
the knot is a torus knot or $K(-1/2,1/3,1/t)$ for odd $t\ge 7$, 
as stated in the last of the proof of Lemma \ref{Lem:LowerBound:RemainingCase}. 
However, in the latter case, 
the surface satisfies $-\chi/\sharp s = 1-2/(t-3) $ and $\sharp s = t-3$. 
Therefore non-torus Montesinos knot have no essential planar surfaces. 
\end{proof}

\begin{proof}[Proof of Corollary \ref{Cor:Denom:UpperBound:Main:ByNonOrientableGenus}]
We only have to check for an essential surface with non-orientable genus $h=1$ and $\sharp b\ge 2$.
In this case, since $\Denom\le(-2+h+\sharp b)/\sharp b+1=2-1/\sharp b<2$,
we have $\Denom=1\le h/2+1=3/2$.

\end{proof}

\subsubsection{The best possibility}

The upper bounds in Theorem \ref{Thm:Denom:UpperBound:Main}, Corollary \ref{Cor:Denom:UpperBound:Main:ByGenus} and Corollary \ref{Cor:Denom:UpperBound:Main:ByNonOrientableGenus} are best possible in a sense.
Let $F$ denote a surface.
Note that we do not care the orientability of the surface $F$.
In this part, $g$ denotes $(2-\chi-\sharp b)/2$,
which coincides with genus if the surface is orientable
and with $2h$ where $h$ is non-orientable genus if the surface is non-orientable.

First, 
we assume that the candidate surface $F$ corresponds to the edgepath system in the item \ref{Item:0-1/2,1-1/2-1/3,B}
in the previous subsection
for odd $t$.
Since the edgepath $\Edgepath_1$ is constant,
$F$ is incompressible by the Proposition 2.1 in \cite{HO},
and thus, is an essential surface.
$-\chi/\sharp s=(\Denom-1)/\Denom$ and $\sharp b=2$ hold,
and give $\Denom=-\chi/\sharp b +1$ and $\Denom=g+1$. 
%
Since $\Denom=(t-3)/2$,
the value of $g=(t-5)/2$ for $t=5,7,9,\cdots$ is $0,1,2,\cdots$.
This indicates that
if $\sharp b\ge2$ are satisfied,
$\Denom\le -\chi/\sharp b +1$ and $\Denom\le g+1$ are best possible for arbitrary non-negative integer $g$.

Next,
let $F$ be a candidate surface in the item \ref{Item:0-1/2,0-1/2-2/3,A}
for odd $t$.
It is incompressible by the Proposition 2.6 in \cite{HO}.
Since $|\Twist|
=4+2/t$, we have $\Denom=t$.
$\sharp s=\mathrm{lcm}(t,t,t)=t$ gives $\sharp b=\sharp s/\Denom=1$.
Thus $-\chi/\sharp s=1$ gives $\Denom=-\chi/\sharp b$ and $\Denom=2g-1$.
%
The value of $g=(t+1)/2$ for $t=3,5,7,\cdots$ is $2,3,4,\cdots$.
This indicates that
if $\sharp b=1$ is satisfied,
$\Denom\le -\chi/\sharp b+1$ and $\Denom\le 2g-1$ are best possible for arbitrary integer $g\ge 2$.
%

%
%

\section{A bound on the difference}
\label{Sec:ABoundOnTheDifference}

The purpose of this section is to prove Theorem \ref{Thm:Diff:UpperBound:Main},
which claims an upper bound of the difference of two boundary slopes.
A large part of this section is the proof of a technical lemma,
which is used for proving the theorem.
The best possibility and some corollaries are also given.

We begin with several remarks on these results. 

\begin{enumerate}

\item
From the argument, 
we must exclude the meridional boundary slope, 
for it corresponds to ``infinity" numerical boundary slope. 
Note that it can actually appear if $\NumTangles\ge4$. 
See \cite{O} for a detail. 

\item
There is an apparent lower bound $|\Slope_1-\Slope_2|\ge 0$.
The lower and upper bounds meet at 
$|\Slope_1-\Slope_2|=0$ and $(-\chi_1/\sharp s_1)+(-\chi_2/\sharp s_2)=-2$. 
This corresponds to the boundary slope of the incompressible disk in 
the trivial knot exterior.

\item
The upper bound  (\ref{Eq:Diff:UpperBound:Main}) is sharp: 
there is an infinite sequence of Montesinos knots 
each of whose exterior includes two essential surfaces with 
boundary slopes satisfying the equality. 
See Subsection \ref{Subsec:Diff:ProofEtc} for more detail.

\item
No such a ``linear" upper bound can hold for $\Delta$. 
See Subsection \ref{SubSec:LinearBoundOfDist} for example. 
In fact known bounds on $\Delta$ are quadratic with respect to 
$-\chi_i/\sharp s_i$. 

\item
In \cite{I}, the first author obtains the same upper bounds for 
2-bridge knot exterior and Seifert fibered manifolds 
which include torus knot exteriors. 
Therefore, the upper bound (\ref{Eq:Diff:UpperBound:Main}) may be 
applicable for some wider class of knot exteriors or manifolds.

\end{enumerate}



\subsection{Linear bound of the distance}
\label{SubSec:LinearBoundOfDist}

Any linear bounds of the distance of two boundary slopes are impossible.
This is a reason why we consider an upper bound of the difference
rather than of the distance.

We give a concrete example of a pair of boundary slopes,
which make any linear bounds impossible.
The example is two boundary slopes $\Slope_1$ and $\Slope_2$ of the Montesinos knot $K(-1/2,1/3,1/t)$ for odd $t\ge 7$.

$\Slope_1$ is a boundary slope which appears in
\ref{Item:0-1/2,1-1/2-1/3,B}
in Subsection \ref{SubSec:ALowerBoundOf-XoverS}.
It corresponds to a type I edgepath system $\EdgepathSystem_1$
\[
\left\{
\begin{array}{l}
\Edgepath_{1,1}=((t-3)/(t-2))\cdot \angleb{-1/2}+(1/(t-2))\cdot \circleb{-1/2} \\
\Edgepath_{1,2}=(((t-5)/(t-3))\cdot \angleb{1/2}+(2/(t-3))\cdot \angleb{1/3})\textrm{\,--\,}\angleb{1/3} \\
\Edgepath_{1,3}=(((t-4)/(t-3))\cdot \angleb{0}+(1/(t-3))\cdot \angleb{1/t})\textrm{\,--\,}\angleb{1/t}.
\end{array}
\right.
\]
For the edgepath system,
$\Twist_1=2/(t-3)$,
$-\chi_1/\sharp s_1=1-2/(t-3)$,
$\Denom_1=(t-3)/2$,
$-\chi_1/\sharp b_1=(t-3)/2-1$.

On the other hand,
the second slope $\Slope_2$ corresponds to a type III edgepath system $\EdgepathSystem_2$
\[
\left\{
\begin{array}{l}
\Edgepath_{2,1}=\angleb{\infty}\textrm{\,--\,}\angleb{0}\textrm{\,--\,}\angleb{-1/2} \\
\Edgepath_{2,2}=\angleb{\infty}\textrm{\,--\,}\angleb{1}\textrm{\,--\,}\angleb{1/2}\textrm{\,--\,}\angleb{1/3} \\
\Edgepath_{2,3}=\angleb{\infty}\textrm{\,--\,}\angleb{1}\textrm{\,--\,}\angleb{1/2}\textrm{\,--\,}\cdots\textrm{\,--\,}\angleb{1/(t-1)}\textrm{\,--\,}\angleb{1/t}.
\end{array}
\right.
\]
For the edgepath system,
$\Twist_2=-2(t+2)$,
$\Denom_2=1$,
$-\chi_2/\sharp s_2=-\chi_2/\sharp b_2=t+2$.

Since
we have both
$(-\chi_1/\sharp b_1)+(-\chi_2/\sharp b_2)=(t-3)/2-1+t+2=(3t-1)/2$
and
$\Delta(\Slope_1,\Slope_2)=\Denom_1 \Denom_2 |\Slope_1-\Slope_2|=\Denom_1 \Denom_2 |\Twist_1-\Twist_2| =t^2-t-5$
at the same time,
$\Delta$ cannot be bounded by an inequality
\[
\Delta(\Slope_1,\Slope_2)\le X(\frac{-\chi_1}{\sharp b_1}+\frac{-\chi_2}{\sharp b_2})+Y
\]
for any constant $X$ and $Y$.


\subsection{An upper bound of the sum of remainder terms}

In this subsection,
we state and prove Lemma \ref{Lem:Remainder:UpperBound},
which claims an upper bound of the sum of ``remainder terms'' of two boundary slopes
and is the key to proving Theorem \ref{Thm:Diff:UpperBound:Main}.


By definition, 
the twist $\Twist$ is roughly twice of a sum of signed lengths of edges.
On the other hand,
as we see in Section \ref{Sec:Preparation},
the major part of $-\chi/\sharp s$ is the sum of lengths of the edgepaths in the edgepath system.
Hence, by the triangle inequality,
these facts imply a bound of the twist
by an inequality like $|\Twist|\le 2\,(-\chi/\sharp s)+\alpha$.
Then,
we introduce the {\em remainder term}
$\Remainder(F)
\equiv
|\Twist|-2\,(-\chi/\sharp s)
$. 
With the remainder term, the key lemma is described as follows.

\begin{lemma}
\label{Lem:Remainder:UpperBound}
For a Montesinos knot $K$,
after simplification,
distinct two candidate surfaces $F_1$ and $F_2$ satisfy
\[
\Remainder_1+\Remainder_2\le 4,
\]
where $\Remainder_i=\Remainder(F_i)$.

In fact, the set of candidate surfaces satisfies following conditions.
\begin{itemize}
\item
For a type I surface $F$, we have $\Remainder\le 4$.
Furthermore, 
there is at most one type I surface with $0<\Remainder\le 4$,
and any other type I surface $F$ satisfies $\Remainder\le 0$.
\item
For a type II surface $F$, we have $\Remainder\le 4$.
Furthermore,
there is at most one type II surface with $0<\Remainder\le 4$,
and any other type II surface $F$ satisfies $\Remainder\le 0$.
\item
For a type III surface $F$, we have $\Remainder\le 0$.
\item
There is at most one surface with $0<\Remainder\le 4$.
Namely, the type I surface $F_1$ with $\Remainder_1>0$ and 
 the type II surface $F_2$ with $\Remainder_2>0$ 
 do not exist for a Montesinos knot $K$ at the same time.
\end{itemize}
\end{lemma}

The simplification is mentioned in Section \ref{Sec:Preparation}.
We divide Lemma \ref{Lem:Remainder:UpperBound} into some partial claims,
and prove the lemma in the rest of this subsection.
We first introduce two notions 
which are used in the proof of Lemma \ref{Lem:Remainder:UpperBound}. 


\subsection*{Cancellation}

On summation in (\ref{Eq:Formula:Twist}),
opposite signs of $\sigma(\Edge)$'s for two or more edges cause cancellation.
If such a cancellation occurs,
we call an edgepath system $\EdgepathSystem$ an {\em edgepath system with cancellation},
and the corresponding surface a {\em surface with cancellation}.
%
A constant edgepath does not cause cancellation.
Only a surface $F$ without cancellation will be able to have $\Remainder>0$ in Lemma \ref{Lem:Remainder:UpperBound}. 

For an edgepath system,
we collect all non-$\infty$-edges of every non-constant edgepath,
divide them into two classes
according to the sign $\sigma(\Edge_{i,j})$ of an edge $\Edge_{i,j}$,
and then sum up the lengths of edges for each class.
With their total lengths $l_{+}$ and $l_{-}$,
let $\Cancel(F)$ denote $\min (l_{+}, l_{-})$.
$\Cancel$ means the amount of the cancellation in calculating the twist for a surface.
With the variable $\Cancel$, the twist is related to the total length of edgepaths in $\Strip$ 
as
\begin{eqnarray}
 |\Twist|
&=& 2\,|l_{+}-l_{-}| = 2\,\{l_{+}+l_{-} -2 \min(l_{+},l_{-}) \}
\nonumber
\\
&=&
2
 \left(
 \sum_{\Edgepath_i \in \NonconstantEdgepaths}
  |\Edgepath_{i,\ge 0} |
 \right)
 -4 \cdot \Cancel
.
\nonumber
\end{eqnarray}


\subsection*{Monotonic edgepath systems}

For an edgepath system, 
the edgepath system is said to be {\em monotonically increasing} (resp. {\em decreasing})
if the $v$-coordinates are monotonically increasing (resp. decreasing) for all edgepaths in the edgepath system.
A surface without cancellation corresponds to a monotonic edgepath system.

Since each vertex $\angleb{p/q}$ ($q>1$) of the diagram $\Diagram$ has two leftward edges, one is increasing and the other is decreasing, there exist only one monotonically increasing basic edgepath and one monotonically decreasing basic edgepath.


\subsubsection{Type II and type III surfaces}



The situation is simplest for type III surfaces among all types of surfaces. 
Even for type II surfaces, the argument is not so complicated.

\begin{lemma}
After simplification,
for any type III surface,
the inequality $\Remainder\le 0$ holds.
\end{lemma}
\begin{proof}
(\ref{Eq:Formula:Twist}) and (\ref{Eq:Formula:EulerCharTypeIII}) give
$
|\Twist|
\le
2\cdot\sum_{i=1}^{\NumTangles}
|\Edgepath_{i,>0}|
=2\,(-\chi/\sharp s)
.
$
Hence, $\Remainder=|\Twist|-2\,(-\chi/\sharp s) \le 0$.
\end{proof}



\begin{lemma}
After simplification,
for any type II surface $F$,
the inequality $\Remainder\le 4$ holds.
Moreover, there exists at most one surface without cancellation, for which $0<\Remainder\le 4$,
while any surface with cancellation satisfies $\Remainder\le 0$.
\end{lemma}
\begin{proof}
(\ref{Eq:Formula:Twist}) and (\ref{Eq:Formula:EulerCharTypeII}) give
\begin{eqnarray}
|\Twist|&\le &
2(\sum_{i=1}^{\NumTangles}|\Edgepath_{i,>0} | )
 +2\,|\EdgepathSystem(+0) |
=2(-\chi/\sharp s+2)
.
 \label{Eq:Twist:UpperBound:TypeII}
\end{eqnarray}
Hence, 
\begin{eqnarray}
\Remainder&=&|\Twist|-2\,(-\chi)/\sharp s \le 4.
\label{Eq:Remainder:UpperBound:TypeII}
\end{eqnarray}

When cancellation occurs,
the difference between the both sides of the inequality (\ref{Eq:Twist:UpperBound:TypeII})
increases by $2\,(+1-(-1))=4$ at a pair of complete edges causing cancellation.
Thus,
$
\Remainder=|\Twist|-2\,(-\chi/\sharp s)\le 0
$.

If the equality in (\ref{Eq:Remainder:UpperBound:TypeII}) holds for an edgepath system,
the edgepath system satisfies either of:
\begin{itemize}
\item
the corresponding basic edgepath system is monotonically decreasing and $\EdgepathSystem(+0)\ge 0$,
\item
the corresponding basic edgepath system is monotonically increasing and $\EdgepathSystem(+0)\le 0$.
\end{itemize}
Both types of edgepath systems are not obtained simultaneously for a Montesinos knot $K$.
By the uniqueness of the monotonically increasing or decreasing basic edgepath system, there is at most one type II surface with $0<\Remainder\le 4$.
\end{proof}


\subsubsection{Type I surfaces}

For type I surfaces, the argument is more complicated than for type II and type III surfaces.
Thus, we here introduce two inequalities for type I surfaces.


\subsection*{An inequality for type I surfaces}

\begin{lemma}
For a type I surface $F$,
its remainder term $\Remainder$ is upper-bounded as
\begin{eqnarray}
\Remainder
&\le&
2(\NumTangles-\NumTangles_{\mathrm{const}})-\left( \NumTangles-2-\sum_{\Edgepath_i \in \ConstantEdgepaths}\frac{1}{q_i} \right)\frac{2}{1-u}
.
\label{Eq:Remainder:UpperBound:TypeI:General}
\end{eqnarray}
\end{lemma}
\begin{proof}
First, an edgepath system for type I surface does not include any vertical edge or $\infty$-edge. 
(\ref{Eq:Formula:Twist}) gives
\begin{eqnarray*}
\Twist&=& 
2 \sum_{i=1}^{\NumTangles}
 \left(
  \left\{
   \begin{array}{l}
    0 \\
    ~~~~~(\textrm{ if $\Edgepath_i$ is constant } ) \\
    \sum_{\Edge_{i,j} \in \Edgepath_{i}} 
     -~\sigma(\Edge_{i,j})~|\Edge_{i,j}| \\
    ~~~~~(\textrm{ if $\Edgepath_i$ is non-constant })
   \end{array}
  \right.
 \right)
.
\end{eqnarray*}

By (\ref{Eq:Formula:EulerCharTypeI}),
\begin{eqnarray}
\Remainder
&=&
|\Twist|-2\frac{-\chi}{\sharp s}
\label{Eq:Remainder:UpperBound:TypeI:WithCancel}
\\
&=&
 -4 \cdot \Cancel
 +2(\NumTangles-\NumTangles_{\mathrm{const}})
 -\left(
 \NumTangles-2-\sum_{\Edgepath_i \in \ConstantEdgepaths}\frac{1}{q_i}
 \right)
 \frac{2}{1-u}
.
\nonumber
\end{eqnarray}
Even if we ignore the effect of cancellation $\Cancel$, 
we have the upper bound in the statement.
\end{proof}


\subsection*{An inequality for type I surfaces with cancellation}

Under the assumption that we could prove $\Remainder\le 4$,
we think about the effect of cancellation.
If $\Cancel \ge 1$,
since cancellation works on the twist by $-4\Cancel$,
$\Remainder\le 0$ immediately follows.
The case of $\Cancel < 1$ only remains. 
In the formula (\ref{Eq:Remainder:UpperBound:TypeI:WithCancel}),
$\Cancel$ and the term $2/(1-u)$ are in the trade-off relationship.
We examine the variation of $\Remainder$ in detail
and make an inequality about $\Remainder$ for an edgepath system with cancellation.

\begin{lemma}
Let $F$ be a type I surface $F$ with cancellation.
Assume that a partial edge of an edge $\angleb{p/q}$\,--\,$\angleb{r/s}$ ($|ps-qr|=1$, $q<s$) causes cancellation.
Then, the remainder term $\Remainder$ is upper-bounded as
\begin{eqnarray}
\Remainder
&\le&
 2(\NumTangles-\NumTangles_{\mathrm{const}}) -2x \min \{s,q+2/x\}
\label{Eq:Remainder:UpperBound:TypeI:WithSingleCancel}
,
\end{eqnarray}
where $x=\NumTangles-2-\sum_{\Edgepath_i \in \ConstantEdgepaths} (1/q_i)>0$.
\end{lemma}
\begin{proof}
Assume that the length of the partial edge included in $\angleb{p/q}$\,--\,$\angleb{r/s}$ is $k/m$, 
and the edge causes cancellation.
We start from (\ref{Eq:Remainder:UpperBound:TypeI:WithCancel}), that is,
$\Remainder
=
 -4 \cdot \Cancel
 +2(\NumTangles-\NumTangles_{\mathrm{const}})
 -x \cdot 2/(1-u)
\le
 -4\, k/m
 +2(\NumTangles-\NumTangles_{\mathrm{const}})
 -x \cdot 2/(1-u)
$%
.
Since the length $k/m$ is given by
$k/m=\{1+s(u-1)\}/\{(s-q)(u-1)\}$ as the formula (\ref{Eq:Formula:LengthOfPartialEdge}),  
the above inequality can be deformed into
$
\Remainder
\le
 2(\NumTangles-\NumTangles_{\mathrm{const}})
 -4 s/(s-q)
 +\{2/(s-q)-x\}\cdot 2/(1-u)
$.
%
The right-hand side is monotonically increasing, constant, or monotonically decreasing as a function of $u$,
according to the sign of $(2/(s-q)-x)$.
Thus, $\Remainder$ can be upper-bounded as in the statement.
%
\end{proof}

Now we show the following.
\begin{lemma}
After the simplification,
for all type I surfaces $F$,
the inequality $\Remainder\le 4$ holds.
Moreover, there exists at most one surface without cancellation, for which $0< \Remainder\le 4$,
while any surface with cancellation satisfies $\Remainder\le 0$.
\end{lemma}

\begin{proof}

First, we here examine type I surfaces without cancellation.
For such a surface, the edgepath system is monotonically increasing or decreasing.
According to the sign of the sum of the tangles $K_i$ as fractions,
only one of the two monotonic basic edgepath systems has a solution of the equation (\ref{Eq:EquationForEdgepathSystem}), that is, $\sum_{i=1}^{\NumTangles} \widetilde{\BasicEdgepath}(u)=0$.
Thus, there exists at most one type I surface without cancellation.

For a type I surface,
since the situation is complicated,
we separate the arguments for $\NumTangles=3$ and $\NumTangles\ge4$.
Furthermore, if $\NumTangles=3$,
we check the lemma case by case 
according to the number of the constant edgepaths.


\medskip
\noindent{(1)
$\NumTangles\ge4$.
}

%

%
First,
without considering the effect of cancellation,
the inequality (\ref{Eq:Remainder:UpperBound:TypeI:General}) is
$\Remainder
\le
2(\NumTangles-\NumTangles_{\mathrm{const}})-\{ \NumTangles-2-\NumTangles_{\mathrm{const}}(1/2) \}\cdot 2/(1-u)
$.
By watching at $u=0$,
we have
$
\Remainder
\le
4-\NumTangles_{\mathrm{const}}
\le
4
$.
%
%
For a surface $F$ with cancellation,
even if $\Cancel<1$,
by (\ref{Eq:Remainder:UpperBound:TypeI:WithSingleCancel}) for $x=\NumTangles-2-\NumTangles_{\mathrm{const}}\cdot 1/2$, 
$q\ge1$ and $s\ge2$,
we obtain
$
\Remainder
\le
\max\{
-2\NumTangles+8,-\NumTangles_{\mathrm{const}}
\}
\le
0
$.

\medskip


\noindent {(2)
$\NumTangles=3$.
}

First, we introduce a notation.
We represent the complete non-horizontal edges including a partial edge $\Edge_i$ of $\Edgepath_i$ by $\angleb{p_i/q_i}$\,--\,$\angleb{r_i/s_i}$,
and the denominators of the $v$-coordinates of constant edgepaths $\Edgepath_a$ and $\Edgepath_b$ by $q_a$ and $q_b$.
Though the indices of non-constant edgepaths may not be successive and be something like $\Edgepath_{1}$ and $\Edgepath_{3}$,
we replace the indices so that the non-constant edgepaths have successive indices like $\Edgepath_{1}$ and $\Edgepath_{2}$,
for ease in the argument.
If cancellation occurs,
$\Edgepath_1$ denotes an edgepath whose partial edge causes cancellation.

\medskip


\noindent {(2-1)
$\NumTangles=3$, $\NumTangles_{\mathrm{const}}=0$.
}

%
%
For an edgepath system with no constant edgepath, the formula (\ref{Eq:Remainder:UpperBound:TypeI:General})
is
$\Remainder\le 6-2/(1-u)$.
Since $0\le u <1$, we have
$\Remainder\le 4$.
%
This is sufficient for edgepath systems without cancellation.
%
%
For an edgepath system with cancellation,
the inequality (\ref{Eq:Remainder:UpperBound:TypeI:WithSingleCancel}) is simplified into
$\Remainder\le 6-2 \min \{s_1,q_1+2\}$.
If $s_1=2$, then $\Remainder\le2$. Otherwise, $\Remainder\le 0$. 
In the case of $s_1=2$, there are two possibilities:
(2-1-1) final edges of all edgepaths have common sign,
(2-1-2) final edges of edgepaths have both positive and negative sign.
In the former case,
except the case $s_2=s_3=2$,
by applying the edgepath system to an inequality similar to (\ref{Eq:Remainder:UpperBound:TypeI:WithSingleCancel}),
we have $\Remainder\le 0$.
For the case $s_1=s_2=s_3=2$,
by solving the equation $\sum_{i=1}^{\NumTangles} \widetilde{\BasicEdgepath_i}(u)=0$,
$\Cancel$ is $1$ or greater, and thus $\Remainder\le 0$.
In the latter case,
the solutions of the equation are non-isolated
and their representative is an edgepath system
which causes no cancellation on $\angleb{p/q}$\,--\,$\angleb{r/s}$.
Thus, we have verified the claim for the remaining case of $s_1=2$ such that there is no such a surface $F$ with cancellation and the remainder term $0<\Remainder\le2$.

\medskip


\noindent {(2-2)
$\NumTangles=3$, $\NumTangles_{\mathrm{const}}=1$.
}

%
%
For a type I surface with one constant edgepath,
the formula (\ref{Eq:Remainder:UpperBound:TypeI:General}) is
$\Remainder
\le
4-(1-1/q_a)\cdot 2/(1-u)
$.
If $q_a=2$, 
since $u \ge 1/2$,
then $\Remainder\le 2$.
If $q_a\ge3$, then $\Remainder\le 0$.
%
%
For an edgepath system with one constant edgepath $q_a=2$ and cancellation,
by the inequality (\ref{Eq:Remainder:UpperBound:TypeI:WithSingleCancel}),
we have $\Remainder\le 4-\min \{s_1,q_1+4\}$,
where the edge
$\angleb{p_1/q_1}$\,--\,$\angleb{r_1/s_1}$ causes cancellation.
The possibility of $\Remainder>0$ remains
when $(q_1,s_1)=(1,2),(1,3),(2,3)$.
In any cases,
we can check that every candidate edgepath system with cancellation
obtained by solving the equation $\sum_{i=1}^{\NumTangles} \widetilde{\BasicEdgepath_i}(u)=0$ actually
satisfies $\Remainder\le 0$ as follows.

First, if the last edge $\Edge_1$ is of type $(q_1,s_1)=(1,2)$,
the $u$ coordinate of the endpoints satisfies $u\le \frac{1}{2}$.
At the same time, as $\Edgepath_a$ is a constant edgepath,
$u\ge \frac{1}{2}$ is required.
Hence, we have a contradiction.

For the case of $(q_1,s_1)=(1,3)$,
possible edgepath systems are essentially same as the following case:
\[
\left\{
\begin{array}{l}
\textrm{$\Edgepath_a$ is a constant edgepath consisting of a point on the edge $\angleb{1/2}$\,--\,$\circleb{1/2}$,} \\
\textrm{$\Edge_1$ is a partial edge of $\angleb{0}$\,--\,$\angleb{1/3}$,} \\
\textrm{$\Edge_2$ is a partial edge of $\angleb{-1}$\,--\,$\angleb{-(t-1)/t}$ ($t=3,4,5$).}
\end{array}
\right.
\]
This appears as the items 
\ref{Item:0-1/2,0-1/3,B},
\ref{Item:0-1/2,0-1/3,C} and
\ref{Item:0-1/2,0-1/3,D}
in the previous section.
Note that $\Edge_1$ and $\Edge_2$ have a common sign.
The value of $|\Edge_1|+|\Edge_2|$ is $1/2+1/2=1$ for $t=3$,
$1/4+1/2=3/4$ for $t=4$.
Since $\Remainder\le 2$ holds even if we ignore the effect of the cancellation, by taking it into account, we have $\Remainder\le 0$ .
When $t=5$, since $|\Edge_1|$ is zero, the edgepath system contradicts the hypothesis that $\Edge_1$ causes cancellation.

For the case of $(q_1,s_1)=(2,3)$,
possible edgepath systems are essentially same as the following case:
\[
\left\{
\begin{array}{l}
\textrm{$\Edgepath_a$ is a constant edgepath consisting of a point on the edge $\angleb{1/2}$\,--\,$\circleb{1/2}$,} \\
\textrm{$\Edge_1$ is a partial edge of $\angleb{1/2}$\,--\,$\angleb{1/3}$,} \\
\textrm{$\Edge_2$ is a partial edge of $\angleb{-1}$\,--\,$\angleb{-(t-1)/t}$ ($t\ge 5$).}
\end{array}
\right.
\]
This corresponds to the item \ref{Item:0-1/2,1-1/2-1/3,B}.
Note that $\Edge_1$ and $\Edge_2$ have opposite signs.
In this case, we have 
$u=(t-1)/\{2(t-2)\}$,
$|\Twist| 
=2/(t-3)$,
and $-\chi/\sharp s 
=1-2/(t-3)$.
Hence, $\Remainder=|\Twist|-2\cdot(-\chi/\sharp s)=6/(t-3)-2$.
If $t=5$, since $|\Edge_1|=0$, $\Edge_1$ does not cause cancellation and the edgepath system contradicts the hypothesis.
If $t\ge 6$, then $\Remainder\le 0$.

\medskip


\noindent {(2-3)
$\NumTangles=3$, $\NumTangles_{\mathrm{const}}=2$.
}

%
%
First, the inequality (\ref{Eq:Remainder:UpperBound:TypeI:General}) is
$\Remainder\le 2-\{ 1-1/q_a-1/q_b \}\cdot 2/(1-u)
$.
Let $q_a$ denote the smaller of denominators of the tangles corresponding to constant edgepaths,
and $q_b$ the larger. 
By the condition for the Montesinos link to be a knot,
we have $q_a\ge 2$ and $q_b\ge3$.
By $u \ge (q_b-1)/q_b$,
we also have $1/(1-u)\ge q_b \ge 3$.
Thus
$\Remainder\le 2-( 1-1/2-1/3 )~2\cdot 3 =1$.
If $q_a=2$ and $q_b=5$, then
$\Remainder\le 2-( 1-1/2-1/5 )~2\cdot 5 = -1$.
If $q_a=3$ and $q_b=3$, then
$\Remainder\le 2-( 1-1/3-1/3 )~2\cdot 3 = 0$.
Similarly, for other edgepath systems, we have $\Remainder\le 0$.

%
%
For an edgepath with cancellation,
it is sufficient to check for $q_a=2$ and $q_b=3$.
Since $x=1-1/q_a-1/q_b=1/6$, 
by inequality (\ref{Eq:Remainder:UpperBound:TypeI:WithSingleCancel}),
$\Remainder\le 2-1/3 \min\{s_1,q_1+12\}$.
If $s_1\ge 6$, then $\Remainder\le 0$.
Hence, an edgepath system with $s_1\le 5$ only remains.
Not so many such concrete examples exist, in fact.
The edgepath system must be
\[
\left\{
\begin{array}{l}
\textrm{$\Edgepath_a$ is a constant edgepath consisting of a point on the edge $\angleb{x/2}$\,--\,$\circleb{x/2}$}, \\
\textrm{$\Edgepath_b$ is a constant edgepath consisting of a point on the edge $\angleb{y/3}$\,--\,$\circleb{y/3}$}, \\ 
\textrm{$\Edge_1$ is a partial edge of $\angleb{p_1/q_1}\textrm{\,--\,}\angleb{r_1/s_1}$} \\
\end{array}
\right.
\]
for some appropriate integer $x$ and $y$.
In order for the edgepath system
to satisfy $\sum_{i=1}^{\NumTangles} \Edgepath_i(u)=0$ at the common $u$-coordinate of the endpoints,
the edge $\Edge_1$ must intersect with the horizontal segments
$v=\pm 1/6+z~(z\in\mathbb{Z})$
within a strip region $2/3 \le u < 1$. 
%
%
The only example of such an edgepath system 
has the partial edge of $\angleb{0}$\,--\,$\angleb{1/5}$ as $\Edge_1$
(or another example essentially same as this example).
We must add at least one increasing complete edge to $\Edge_1$
so that $\Edge_1$ actually causes cancellation.
Since $|\Edge_1|=1/2$,
we have $\Remainder=1-4\cdot 1/2=-1$ for this edgepath system.

\medskip


\noindent {(2-4)
$\NumTangles=3$, $\NumTangles_{\mathrm{const}}=3$.
}

As mentioned in Remark \ref{Rem:Condition:ForLinkToBeKnot},
no edgepath system with three constant edgepaths exists.
\end{proof}


\subsubsection{Type I surface and type II surface without cancellation}

Now, we have only to show the following.

\begin{lemma}
The type I surface without cancellation and
the type II surface without cancellation 
do not exist for a Montesinos knot $K$ at the same time.
\end{lemma}
\begin{proof}
For the type I surface without cancellation,
$\sum_{i=1}^{\NumTangles} K_i$ and $\EdgepathSystem(+0)$ have opposite signs.
In contrast with this,
for the type II surface without cancellation,
$\sum_{i=1}^{\NumTangles} K_i$ and $\EdgepathSystem(+0)$ have the same sign or $\EdgepathSystem(+0)=0$.
\end{proof}

Combining sub-lemmas completes the proof of Lemma \ref{Lem:Remainder:UpperBound}.
%


\subsection{Proof, best possibility and corollaries of Theorem \ref{Thm:Diff:UpperBound:Main}}
\label{Subsec:Diff:ProofEtc}

\subsubsection{Proof}

Theorem \ref{Thm:Diff:UpperBound:Main} follows Lemma \ref{Lem:Remainder:UpperBound}.
\begin{proof}[Proof of Theorem \ref{Thm:Diff:UpperBound:Main}]
By the triangle inequality and Lemma \ref{Lem:Remainder:UpperBound}
\begin{eqnarray*}
|\Slope_1-\Slope_2| 
&=& |\{\Twist_1-\Twist_{S}\} -\{\Twist_2-\Twist_{S}\} | \\
&=& |\Twist_1-\Twist_2| \\
&\le& |\Twist_1|+|\Twist_2| \\
&=&2\,(\frac{-\chi_1}{\sharp s_1}+\frac{-\chi_2}{\sharp s_2})+(\Remainder_1+\Remainder_2) \\
&\le&2\,(\frac{-\chi_1}{\sharp s_1}+\frac{-\chi_2}{\sharp s_2})+4
,
\end{eqnarray*}
where $\Twist_{S}$ is the twist of a Seifert surface $F_{S}$ of $K$.

Note that there is no side effect by simplification.
For example, the inequality (\ref{Eq:Diff:UpperBound:Main})
holds also for two surfaces $F_1$ and $F_2$ with the same boundary slope,
since essential surfaces satisfy $(-\chi_1/\sharp s_1)+(-\chi_2/\sharp s_2)\ge -2$ even if they are disks.

.
\end{proof}

\subsubsection{Best possibility}

The upper bound in Theorem \ref{Thm:Diff:UpperBound:Main} is best possible in a sense.

First, we provide a concrete example of a family of Montesinos knots and pairs of edgepath systems of boundary slopes of the knots. 
The Montesinos knot is $K(1/(2k),1/5,\ldots,1/5)$ with $\NumTangles\ge 3$ tangles and the natural number $k$.
The two edgepath systems $\EdgepathSystem_1=(\Edgepath_{1,1},\Edgepath_{1,2},\ldots,\Edgepath_{1,\NumTangles})$ and $\EdgepathSystem_2=(\Edgepath_{2,1},\Edgepath_{2,2},\ldots,\Edgepath_{2,\NumTangles})$ of $F_1$ and $F_2$ are as follows.
\begin{eqnarray*}
\Edgepath_{1,1}&=&\angleb{0}-\angleb{1/(2k)}, \\
\Edgepath_{1,i}&=&\angleb{0}-\angleb{1/5} \textrm{\ \ \ \ for $2\le i \le \NumTangles$}, \\
\Edgepath_{2,1}&=&\angleb{\infty}-\angleb{1}-\angleb{1/2}-\angleb{1/3}-\cdots-\angleb{1/(2k-1)}-\angleb{1/(2k)}, \\ 
\Edgepath_{2,i}&=&\angleb{\infty}-\angleb{1}-\angleb{1/2}-\angleb{1/3}-\angleb{1/4}-\angleb{1/5} \textrm{\ \ \ \ for $2\le i \le \NumTangles$}. 
\end{eqnarray*}
These candidate edgepath systems $\EdgepathSystem_1$ and $\EdgepathSystem_2$ are type II and type III respectively.
Two candidate surfaces $F_1$ and $F_2$ are both incompressible by the Corollary 2.4 and Proposition 2.5 in \cite{HO}.


Since $\EdgepathSystem_1$ and $\EdgepathSystem_2$ are monotonically decreasing and increasing,
their remainder term $\Remainder$ are easily confirmed to be $4$ and $0$.
$\Twist_1$ and $\Twist_2$ with opposite signs give $|\Twist_1-\Twist_2|=|\Twist_1|+|\Twist_2|$.
Thus, $F_1$ and $F_2$ satisfy
\[
|\Slope_1-\Slope_2|= 2\,(\frac{-\chi_1}{\sharp s_1}+\frac{-\chi_2}{\sharp s_2})+4
.
\]
Note that $(-\chi_1/\sharp s_1)+(-\chi_2/\sharp s_2)$
is greater than arbitrary $t$ for sufficiently large $k$.

\subsubsection{Corollaries}

As described in Section \ref{Sec:ABoundOnTheDenominator},
if a Montesinos knot $K$ is not the same as or is not isotopic to $(-2,3,t)$-pretzel knots for odd $t\ge 3$,
we have $-\chi/\sharp s\ge 1$ for its boundary slopes.
The $(-2,3,t)$-pretzel knots have boundary slopes and corresponding essential surfaces as follows. Note that they are torus knots if $t=3$ or $5$.

\begin{itemize}
 \item $t=3$
 \begin{itemize}
  \item[$\circ$] $\Slope_a=12$, $\Denom_a=1$, $\chi_a=0$, $\sharp s_a=\sharp b_a=2$, $-\chi_a/\sharp s_a=-\chi_a/\sharp b_a=0$,
  \item[$\circ$] $\Slope_b=0$, $\Denom_b=1$, $\chi_b=-5$, $\sharp s_b=\sharp b_b=1$, $-\chi_b/\sharp s_b=-\chi_b/\sharp b_b=5$.
 \end{itemize}
 \item $t=5$
 \begin{itemize}
  \item[$\circ$] $\Slope_a=15$, $\Denom_a=1$, $\chi_a=0$, $\sharp s_a=\sharp b_a=2$, $-\chi_a/\sharp s_a=-\chi_a/\sharp b_a=0$,
  \item[$\circ$] $\Slope_b=0$, $\Denom_b=1$, $\chi_b=-7$, $\sharp s_b=\sharp b_b=1$, $-\chi_b/\sharp s_b=-\chi_b/\sharp b_b=7$.
 \end{itemize}
 \item $t\ge 7$,
 \begin{itemize}
  \item[$\circ$] $\Slope_a=16$, $\Denom_a=1$, $\chi_a=6-t$, $\sharp s_a=\sharp b_a=1$, $-\chi_a/\sharp s_a=-\chi_a/\sharp b_a=t-6$,
  \item[$\circ$] $\Slope_b=(t^2-t-5)/\{(t-3)/2\}$, $\Denom_b=(t-3)/2$, $\chi_b=5-t$, $\sharp s_b=t-3$, $\sharp b_b=2$, $-\chi_b/\sharp s_b=1-2/(t-3)$, $-\chi_b/\sharp b_b=(t-5)/2$,
  \item[$\circ$] $\Slope_c=2t+6$, $\Denom_c=1$, $-\chi_c/\sharp s_c=-\chi_c/\sharp b_c=1$,
  \item[$\circ$] $\Slope_d=0$, $\Denom_d=1$, $-\chi_d/\sharp s_d=-\chi_d/\sharp b_d=t+2$.
 \end{itemize}
\end{itemize}
By examining these boundary slopes,
we obtain linear or quadratic upper bounds of the difference and the distance of the two boundary slopes.




\begin{proof}[Proof of Corollary \ref{Cor:Diff:UpperboundByGenus}]

If the knot $K$ is not $(-2,3,t)$-pretzel,
$-\chi/\sharp s\ge 1$ and $g\ge 1$ hold for all essential surfaces.
Since $g\ge 1$ gives $-\chi/\sharp s\le 2g-1$,
we have the inequality (\ref{Eq:Diff:UpperBound:Main-B:ByGenus}).
If the knot $K$ is $(-2,3,t)$-pretzel with $t\ge 7$,
then genus of any essential surface is found to be $1$ or greater.
Similarly to the previous case, we have (\ref{Eq:Diff:UpperBound:Main-B:ByGenus}).
For the remaining $(-2,3,3)$ and $(-2,3,5)$-pretzel knots,
the value of $|\Slope_1-\Slope_2| - 4 \,( g_1 + g_2 )$ 
for the two boundary slopes 
is $|12-0|-4\,(0+3)=0$ for $(-2,3,3)$
and $|15-0|-4\,(0+4)=-1<0$ for $(-2,3,5)$.
\end{proof}



\begin{proof}[Proof of Corollary \ref{Cor:DiffDist:UpperBound}]
If the knot $K$ is neither $(-2,3,3)$ nor $(-2,3,5)$-pretzel essentially, we have $-\chi/\sharp s\ge 1/2$ for all boundary slopes.
Then, since $(-\chi_1/\sharp s_1)+(-\chi_2/\sharp s_2)\ge 1$,
we obtain 
(\ref{Eq:Diff:UpperBound:Linear})
from (\ref{Eq:Diff:UpperBound:Main}) in Theorem \ref{Thm:Diff:UpperBound:Main}.

For remaining $(-2,3,3)$ and $(-2,3,5)$-pretzel knots,
since the expression $|\Slope_a-\Slope_b|-6\,((-\chi_a/\sharp s_a)+(-\chi_b/\sharp s_b))$
has its value $12-6\cdot 5=-18<0$ and $15-6\cdot 7=-27<0$ respectively,
the inequality (\ref{Eq:Diff:UpperBound:Linear}) holds.
\end{proof}





\begin{proof}
[Proof of Corollary \ref{Cor:Main3b}]
If both boundary slopes satisfy $-\chi/\sharp s\ge 1$,
or equivalently $\Denom \le -\chi/\sharp b$,
we obtain 
(\ref{Eq:Dist:UpperBound:Quadratic})
easily from (\ref{Eq:Dist:UpperBound:Main-C}).
Thus, we are done for all Montesinos knots but $(-2,3,t)$-pretzel knots.

For $(-2,3,3)$ and $(-2,3,5)$-pretzel knots,
since one of the two essential surfaces has Euler characteristic $0$,
there are no pairs of boundary slopes to be applied to the inequality.
For $(-2,3,t)$-pretzel knots with odd $t\ge 7$,
the value of the expression $\Delta(\Slope_i,\Slope_j)-8\cdot (-\chi_i/\sharp b_i)\cdot (-\chi_j/\sharp b_j)$ for $(i,j)=(a,b),(b,c),(b,d)$ are $-3t^2+35t-101$, $-3t+16$ and $-3t^2+11t+35$ respectively, which are all negative for any $t\ge 7$.
\end{proof}

%
%

\section{Known results and open problems}

Here we give a brief review about 
the study of boundary slopes of essential surfaces 
related to our results, and collects some open problems.

\subsection*{}

For the existence and the number of boundary slopes, 
the following are fundamental. 
It was shown by Hatcher in \cite{H} that 
there are just finitely many boundary slopes 
of essential surfaces for a compact, orientable, irreducible $3$-manifold 
with boundary a single torus. 
Also it was shown by Culler and Shalen in \cite{CS84} that 
there are at least two boundary slopes 
for a non-trivial knot exterior in the $3$-sphere $S^3$. 
See also \cite{CS04}. 

Boundary slopes for some class of knots 
have been intensively studied. 
As a pioneering work, 
for two-bridge knots, Hatcher and Thurston gave a complete enumeration 
of boundary slopes in \cite{HT}. 
Following this work, 
Hatcher and Oertel \cite{HO} developed a procedure to compute 
the boundary slopes for Montesinos knots, 
on which our arguments heavily depend. 
Recently boundary slopes of genus one essential surface 
for Montesinos knots of length three are completely determined 
by Wu \cite{W}.

\subsection*{}

The denominators of boundary slopes have also been studied
in relation to the study of Dehn surgery.

In the following,
let $F$ be an essential surface properly embedded
in the exterior of a non-trivial knot $K$ in $S^3$.
The surface $F$ is of Euler characteristic $\chi$ and
has the boundary slope $\Slope$ of $F$,
which is represented by an irreducible fraction $\Slope = \Numer / \Denom$.
Let $\sharp s$ denote the number of sheets of $F$ and
$\sharp b$ the number of boundary components of $F$,
where they are related to each other by $\sharp s = \Denom \,\sharp b$.

Please note that
the results cited below will often be modified from the original  
statements.
It is for making easy to see their relationship and to compare with our  
results.

First, for genus 0 case, by Gordon and Luecke in \cite{GLu87},
it was shown that $\Denom \le 1$, that is, $\Slope$ is integral.
On the other hand, for genus one case,
no corresponding results are found in literature at least by the author.
The related result was obtained by Gordon and Luecke in \cite{GLu95,  
GLu00}.
That is,
if a Dehn surgery on a hyperbolic knot in $S^3$ along a slope $\Slope$
yields a closed $3$-manifold containing an incompressible torus, then  
$\Denom \leq 2$.
Note that if such a surgery can occur, then $\Slope$ is a boundary slope of
an essential surface of genus one. 
However the converse does not hold in general. 

Following these results, it is natural to ask:

\begin{problem}
Find a generalization to these results 
for the higher genus case or the non-orientable surface case. 
\end{problem}

Concerning this problem, some results are already known. 
By using the argument used in \cite[Proposition 6.1]{GLi},
together with Gabai's thin position argument \cite{Ga87},
we have

$$ \Denom \le 6 \frac{ - \chi }{ \sharp b } $$
if the knot $K$ is non-cabled.
Originally in their argument, $F$ is assumed to be orientable,
but the assumption might be not necessary.
This result had not been included in \cite{GLi}, but was suggested in  
\cite{R}.

In \cite{T}, other generalization was developed
if $K$ admits some tangle decomposition.
In particular if $K$ has non-trivial $t$ connected summands, he obtained

$$ \Denom \le \frac{ g }{ t-1 } ,$$
where $F$ is assumed to be orientable and $g$ denotes the genus of $F$.

Under restriction to the class of knots, two excellent results are  
known.
One is for two-bridge knots, by Hatcher and Thurston \cite{HT}.
They gave a classification of essential surfaces in two-bridge knot  
exteriors, and
as a corollary, it was shown that all such surfaces have integral  
slopes.
Another one is for alternating knots, by Menasco and Thistlethwaite  
\cite{MT}.
They presented that

$$ \Denom \le \frac{ - \chi }{ \sharp b } $$
for non-torus alternating knots.
As a corollary they achieved the affirmative answer to
the well-known Cabling Conjecture for alternating knots.
We also remarked that
it is known that a torus knot exterior contains only two essential  
surfaces
and their boundary slopes are both integers.

Recently, Matignon and Sayari \cite{MS04}obtained
similar bounds for non-orientable surfaces.
Their result was obtained by using Dehn surgery method, 
and in fact, they do not assume 
that the surfaces they are considering are essential. 
However, as pointing out in \cite{Te}, 
the condition that the surfaces are essential is necessary. 
Under the essentiality condition, 
their results could be interpreted in terms of boundary  slopes as follows. 
They actually showed that

$$ \Denom \le \frac{ - \chi }{ \sharp b } + 4 $$
if $F$ is non-orientable and $\sharp b > 1$.
If $\sharp b = 1$, they had

$$ \Denom \le  - 3 \chi  + 1 $$
if $K$ is not a cable knot, and

$$ \Denom \le  - 5 \chi  + 3 $$
if $K$ is cable knot.
They also showed that

$$ \Denom \le \frac{-\chi}{\sharp b} $$
if $K$ is a composite knot,
and

$$ \Denom \le \frac{-\chi}{\sharp b} +1$$
if $K$ admits a Conway sphere.

Remark that, for a non-orientable genus two case, namely,
punctured Klein bottle case,
it was obtained in \cite{GLu95} that $\Denom = 1$, that is, $\Slope$ is integral.\\

There are many results which give upper bounds on 
the distances between boundary slopes. 
The main problem would be:

\begin{problem}
Establish a sharp estimate on the distances between boundary slopes 
in terms of the genera of the corresponding essential surfaces. 
\end{problem}

In the following let $M$ be a compact orientable irreducible $3$-manifold
whose boundary $\partial M$ is homeomorphic to the torus $T^2$.
For $i =1,2$, let $F_i$ be an essential surface properly embedded in $M$ 
of Euler characteristic $\chi_i$. 
The boundary slope $\Slope_i$ of $F_i$ is 
represented by an irreducible fraction $\Slope_i=\Numer_i/\Denom_i$.
Let $\sharp s_i$ denote the number of sheets.
Note that 
if the number of boundary of $F_i$ is denoted by $\sharp b_i$,
$\sharp b_i$ and $\sharp s_i$ are related to each other by $\sharp s_i = \Denom_i \sharp b_i$.

For small genus surface case, intensively fine results have been achieved 
in relation to the study of the exceptional Dehn surgery. 
If both $F_i$'s are planar, 
Gordon and Luecke proved in \cite[Theorem~1.1]{GLu96} that $\Delta \leq 4$ holds. 
If both $F_i$'s are punctured torus, Gordon proved in \cite[Theorem~1.1]{G} 
that $\Delta \leq 8$ holds. 
Moreover, he gave $\Delta \leq 5$ with just five exceptional manifolds, 
which are completely characterized. 

As a generalization to the higher genera case, 
Gordon and Litherland obtained 
in \cite[Proposition~6.1]{GLi} the following: 
Suppose that $M$ contains no cable spaces. 
If both $F_i$'s are orientable and $F_1$ is planar, 
then $ \Delta < 6 \,(-\chi_2/\sharp b_2) $ holds. 

As a natural extension of \cite[Proposition~6.1]{GLi}, 
Torisu obtained in \cite[Theorem~1]{T} the following: 
Suppose that $M$ contains no essential annuli. 
If both $F_i$'s are orientable and of genus $g_i \geq 1$, 
then $ \Delta < 36 (2 g_1 -1 ) (2 g_2 -1)$ holds. 

On the other hand, Rieck obtained in \cite{R} a slightly sharper bound. 
Suppose that $M$ contains no essential annuli. 
If both $F_i$'s are orientable and of genus $g_i \geq 1$, 
$ \Delta < 18 (2 g_1 +1 ) (g_2 +1)$ holds. 
Moreover if $\sharp b_i \geq 2$ for $i=1,2$, 
then $ \Delta < 18 (g_1 +1 ) (g_2 +1)$ holds. 
This is slightly different from the original form. 
Please refer \cite[Theorem~5.2]{R} as the original form. 
In fact, from his proof, we can find 
$$ \Delta < 
18 \left(2 \frac{g_1}{\sharp b_1} +1 \right) 
\left(2 \frac{g_2}{\sharp b_2} +1 \right).$$ 

These above are all proved by 
the combinatorial analysis of the graph 
constructed from the intersection of the two surfaces. 

On the other hand, by a differential geometric approach, 
the following bound was shown by Hass, Rubinstein and Wang in \cite{HRW}: 
Suppose that the interior of $M$ admits 
a complete hyperbolic metric of finite volume. 
Then 
$$ \Delta \leq \frac{(2\pi)^2}{3.35} 
\frac{ - \chi_1 }{\sharp b_1} \frac{ - \chi_2 }{\sharp b_2}
\doteq 11.8 
\frac{ - \chi_1 }{\sharp b_1} \frac{ - \chi_2 }{\sharp b_2}
$$
holds. 
This bound still holds in the case that 
the surface $F_i$ is immersed essential surface. 
We remark that this also differs from the original form. 
In the original form, the surfaces are assumed to be orientable, 
but this orientability condition is not necessary in their argument. 
Please refer \cite[Theorem~4.5]{HRW} as the original form. 
Moreover, by using the result of Agol \cite[Theorem~5.1]{A}, 
this bound is improved as 
$$ \Delta \leq \frac{36}{3.35} 
\frac{ - \chi_1 }{\sharp b_1} \frac{ - \chi_2 }{\sharp b_2}
\doteq 
\frac{43}{4} 
\frac{ - \chi_1 }{\sharp b_1} \frac{ - \chi_2 }{\sharp b_2} .
$$

%
%

\end{document}